\documentclass{amsart}

\usepackage{amssymb,amsxtra,amsthm, amsfonts}
\usepackage{amscd}
\usepackage{epsfig}
\input xy
\xyoption{all}

{\end{list}}

\numberwithin{equation}{section} 

\theoremstyle{plain}
\theoremstyle{plain}
\newtheorem*{acknowledgement} {Acknowledgement}
\newtheorem{proposition}{Proposition}[section]

\newtheorem{theorem}[proposition]{Theorem}

\newtheorem{lemma}[proposition]{Lemma}
\newtheorem{corollary}[proposition]{Corollary}

\theoremstyle{definition}
\newtheorem{definition}[proposition]{Definition}
\newtheorem*{claim}{Claim}
\newtheorem*{case}{Case} 

\newtheorem{remark}[proposition]{Remark}
 
\newtheorem{chunk}{Case}
\addtocounter{chunk}{-1}
\newtheorem{chunkk}{}
\numberwithin{chunkk}{subsubsection} 

\newcommand{\A}{{\mathcal A}}
\newcommand{\C}{{\mathbb C}}

\newcommand{\DD}{{\tilde D}} 
\newcommand{\D}{{\Delta}}
 
\newcommand{\Ex}{{\mathcal Ext^1_{\OO_D}}}
\newcommand{\EE}{{\mathcal E}} 
\newcommand{\F}{{\mathcal F}}
\newcommand{\G}{{\Gamma}} 
\newcommand{\I}{{\mathcal I}}
\newcommand{\II}{{\mathcal I_D/\mathcal I_D^2 }} 

\newcommand{\LL}{{\mathcal L}} 
\newcommand{\N}{{\mathcal N}}
\newcommand{\PP}{{\mathbb P}} 
\newcommand{\OO}{{\mathcal O}}
\newcommand{\Om}{{\Omega^1_D}} 
\newcommand{\T}{{\mathcal T}}
\newcommand{\U}{{\mathcal U}} 

\newcommand{\X}{{\mathcal X}}

\newcommand{\lright}{\mathop{\longrightarrow}} 
\newcommand{\rright}{\mathop{\rightarrow}}
\newcommand{\Spec}{\operatorname{Spec}} 
\newcommand{\Hom}{\operatorname{Hom}}
\newcommand{\Der}{\operatorname{Der}} 

\newcommand{\Ext}{\operatorname{Ext}} 
\newcommand{\Exal}{\operatorname{Exalcom}}
\newcommand{\im}{\operatorname{im}} 
\newcommand{\coker}{\operatorname{coker}}
\newcommand{\mult}{\operatorname{mult}} 
 \newcommand{\Bl}{\operatorname{Bl}}

\newcommand{\al}{{\alpha}}
\newcommand{\bb}{{\beta}} 
\newcommand{\dd}{{\delta}}
\newcommand{\f}{{\varphi}} 
\newcommand{\om}{{\omega_D}}
\newcommand{\s}{{\sigma}} 
\newcommand{\codim}{\operatorname{codim}}

\begin{document} 

\title{Analytic rigidity of $K$-trivial extremal 
contractions of smooth 3-folds} 

\author{Csilla Tam\'as}
\thanks{This paper is part of the author's Ph.D. thesis at Purdue
University, written under the guidance of K. Matsuki.} 
\address{Department of Mathematics\\ University of Georgia\\ Athens, GA 30602-7403\\ USA} 
\email{ctamas@math.uga.edu}

\date{\today}

\begin{abstract}
We discuss the problem of classifying birational
extremal
contractions of smooth threefolds where the canonical bundle is trivial
along the curves contracted, in the case when a divisor is contracted to a
point. We prove the analytic rigidity of the contraction in the case this
divisor is normal of degree $\geq 5$ (i.e. we show that the analytic structure
of the contraction is
completely determined by the isomorphism class of the exceptional locus
and its
normal bundle in X). This was previously known only in the case of smooth
exceptional locus.
\end{abstract}

\maketitle


\section{Introduction} 

In the minimal model program, the study of certain types of birational contractions, called {\it extremal},
is of central importance. In \cite{Mori}, S. Mori 
studied and classified birational extremal
contractions $\varphi \colon X \longrightarrow~Y$ of smooth threefolds
$X$ where the canonical bundle of $X$ is {\it negative} along the
fibers of the contraction.
His classification in the case the exceptional locus $Exc(\f)$ is a divisor
includes the following results: When $Exc(\f)$ contracts to a curve on
$Y$, $Exc(\f)$ is a $\PP^1$-bundle over the base curve and $Y$ is smooth.  When  $Exc(\f)$ contracts to a point $q \in Y$, it is either $\PP^2$,
$\PP^1 \times \PP^1$ or a singular quadric, with specified normal
bundle, and $X$ is the blow-up of $q$ on $Y$; in this case the
analytic structure of the neighborhood of $Exc(\f)$ is uniquely
determined by the isomorphism class of $Exc(\f)$ and its normal bundle
in $X$. We call this feature the {\it analytic rigidity} of
the contraction (see Definition \ref{Def:an.rig.}). 

In this paper we attempt to give a similar description for birational
extremal contractions of smooth 3-folds in the {\it $K$-trivial} case,
i.e. where the
canonical bundle is numerically trivial on all the curves
contracted (see \S \ref{known} for a precise definition). Our main result is
\begin{theorem}[Main Theorem]\label{T:main} 
Let $X$ be a smooth projective 3-fold over $\C$  and let $\varphi \colon X 
\longrightarrow~Y$ be a $K$-trivial birational extremal contraction onto a normal 
projective variety $Y$, contracting a divisor $D \subset X$ to a point $q \in Y$. Suppose $D$ is 
normal, and $(\omega_D^2) \geq 5$. Then the contraction $\f$ is 
analytically rigid. 
\end{theorem} 

This way we obtain the classification of $K$-trivial
extremal contractions in terms of the exceptional divisor $D$ and its
normal bundle $\N_{D/X}$ in $X$ in the case when $D$ is normal with
$d=(\omega_D^2) \geq 5$, contracting to a point. As the exceptional locus is a normal
rational (possibly singular) del Pezzo surface $D$ of degree $d \geq 5$
with normal bundle isomorphic to $\om$ (see \S \ref{known}),
we obtain a finite list of possible contractions up to analytic
isomorphism.

Although  extremal contractions are algebraic objects, establishing
algebraic classification would lead to the problem of finding moduli of algebraic objects, which is global in
nature, while the classification of extremal contractions should be
local in nature. From this point of view, we may regard analytic classification as the best
possible result. 

As the result of analytic rigidity shows, in order to know the 
analytic structure of the contraction, it is sufficient to construct one 
example for each possible exceptional locus. No algebraic examples are 
known except when $D$ is a nonsingular del Pezzo surface of degree 6 
(\cite{Nam1}). In \S \ref{examples} we  present an
example of an embedding of a singular del Pezzo surface of degree 7
into a smooth projective 3-fold such that the canonical bundle of the 
ambient space is numerically trivial on $D$. By Fujiki's 
contraction theorem, $D$ can be contracted analytically, giving us an
analytic $K$-trivial contraction.  Similar constructions can be
carried out  for each possible exceptional 
divisor $D$ (any normal rational del Pezzo 
  surface of degree $d \geq 5$).

$K$-trivial extremal contractions have been studied by M. Gross (\cite{Gr}) and P.M.H.~Wilson (\cite{Wilson}). 
Their analysis
in the case when the exceptional locus $D$ is a divisor contracting to a point includes a description
of the (possible) exceptional loci, and the description of the
analytic structure (and in particular establishing the analytic
rigidity) in the case when $D$ is nonsingular with
$(\omega_D^2) \geq 5$. The analytic rigidity in the case when $D$ is non-normal or $(\omega_D^2)
\leq 4$ is not known. Wilson also describes the possible contractions in
the case when the exceptional locus $D$ contracts to a
curve. Small $K$-trivial   birational extremal contractions of smooth 3-folds (i.e.
when the exceptional locus is a collection of curves ) are
three-dimesional flopping contactions; these were  studied
in  \cite{Kollar} and \cite{Reid2}.

Even in the surface case, the condition of  $K$-negativity or $K$-triviality is essential
for analyitic rigidity. In \cite{Laufer}, Laufer gives a complete list of taut
surface singularities (i.e. normal 2-dimensional singularities that
are completely determined by the weighted dual graph $\G$ of the
exceptional locus $E$ of the minimal resolution). He also lists those
singularities which are determined by the weighted dual graph and the
analytic structure of $E$, and states that the singularities obtained
by contractions
of curves of general type are not determined by $\G$ and the
analytic structure of $E$.

We outline below the proof of the Main Theorem:
 
By \cite[Th. 3]{Hir}, the proof of analytic rigidity 
is reduced to showing that any two embeddings of $D$ into smooth 
complex 3-folds with normal bundles isomorphic to $\OO_D(K_D)$ are 
formally equivalent. Suppose now that we have two $K$-trivial extremal
contractions $\f \colon X \lright Y$ and $\f' \colon X' \lright Y'$,
with  isomorphic exceptional divisors $D$ and $D'$, and such that
$\N_{D/X} \simeq \N_{D'/X'}$. To prove the formal equivalence, we
first show that if $H^1(D, \T_D \otimes \II)=0$, then the two
embeddings $D \subset X$ and $D' \subset X'$ are 2-equivalent, i.e. we
have an isomorphism of the ringed spaces $(D, \mathcal  O_X/\mathcal
I_D ^{2})$ and  $(D', \mathcal  O_{X'}/\mathcal I ^2_{D'})$, where $\I_D$ denotes
the ideal sheaf of $D$ in $X$. Then we
can obtain a formal equivalence by showing that the obstruction spaces
$H^1(D, \mathcal T_X \otimes \mathcal I ^{\nu}/\mathcal I 
^{\nu +1})$
to extending a $\nu$-equivalence ($\nu \geq 2$) to a $\nu
+1$-equivalence vanish for any $\nu \geq 2$. The vanishing of both
$H^1(D, \T_D \otimes \II)$ and $H^1(D, \mathcal T_X \otimes \mathcal I_D ^{\nu}/\mathcal I _D
^{\nu +1})$ is reduced to showing $H^1(D, \T_D)=0$ using properties of
del Pezzo surfaces.  This last vanishing is then proved using an
explicit description of normal rational del Pezzo surfaces.

\begin{remark} 
 We should note that our results about formal equivalence hold over any
algebraically closed field of characteristic 0. However, over an arbitrary  
field  there is no notion of analytic rigidity (and in particular we
dont't have \cite[Th. 3]{Hir}).
 Over an
arbitrary field, formal equivalence implies only 
equivalence in the \'etale topology (\cite[Th. (4.6)]{Artin1}). 
\end{remark}

\noindent{\bf Conventions/Notations.}
\begin{enumerate} 
\item We are working over the field of complex numbers $\C$.
\item For  an algebraic variety $D$, we denote by $\T_D$ its tangent
sheaf $\mathcal Hom_{\OO_D}(\Om,\OO_D)$. 
 In case it
exists (for example if $D$ is Cohen-Macaulay), we denote the canonical sheaf of
$D$  by $\om$. If $\om$ is locally free, we denote  by $K_D$ the (linear equivalence class of the) canonical
divisor of $D$. 
\item For a coherent sheaf $\F$ on $D$, we denote $h^i(D,\F) =\dim
H^i(D,\F)$.
\item We denote the intersection number of a (Cartier) divisor $\D$
and a curve $C$ by $(\D \cdot C)$. 
\item For any variety $X$, we denote by $\overline{NE}(X)$ the
closure of the cone of effective 1-cycles on $X$.
\end{enumerate}

\section{Preliminaries}\label{known}  
  
\subsection{$K$-trivial extremal contractions}
\begin{definition} Let $X$ be a smooth projective $n$-fold over
$\mathbb C$ and $\varphi \colon X \longrightarrow~Y$ a
birational morphism onto a normal projective variety $Y$ such that the
exceptional locus $D$ of $\f$  is of codimension 1 and such that $\dim \f (D)=0$.
We call the contraction $\f$ 
\begin{description} 
  \item[(P1)] \underline{extremal} if all the curves
  contracted by $\varphi$ are numerically proportional, i.e. given two curves $C$ and $C'$ contracted by $\f$, 
there is a rational number $r$ such that for any divisor $E$ in $X$, 
we have $(E \cdot C')=r(E \cdot C)$.
  \item[(P2)] \underline{$K$-trivial} if the canonical
  bundle on $X$ is numerically trivial on all curves
  contracted by $\varphi$, i.e. $(K_X \cdot C)=0$ for any curve $C$
  contracted by $\f$ to a point.
\end{description} 
\end{definition}
Note that the condition (P1) implies that the exceptional locus 
  $D$ is an irreducible divisor (Proposition \ref{P:aboutD}), and hence $\f$
contracts $D$ to a point $q \in Y$.
\begin{definition}\label{Def:an.rig.}
The contraction $\f$ is called \underline{analytically rigid} if its
analytic structure is uniquely determined by the isomorphism class
of $Exc(\f)=D$ and its normal bundle $\N_{D/X}$ in $X$. More precisely, suppose 
$\f' \colon X' \lright Y'$ is another birational map on a smooth projective 3-fold $X'$
with exceptional locus $D'$, contracting $D'$ to a point $q' \in Y'$. If 
$D \simeq D'$ and we have an
isomorphism of normal bundles $\N_{D/X} \simeq \N_{D'/X'}$, then
analytic rigidity of $\f$ means that there
are open (analytic) neighborhoods $U$ of $D$ in $X$ and $U'$ of $D'$ in $X'$ over
which the contractions $\f$ and $\f'$ are analytically isomorphic,
i.e. we have the following commutative diagram (in the analytic category):
\[
\xymatrix{D\ar[d] &\subset &U\ar[d] &\simeq &U'\ar[d] &\supset &D'\ar[d]\\
q &\in &\f(U) &\simeq &
\f'(U') &\ni &q' }
\]
\end{definition}

\begin{remark}
A priori our definition (P1) of an extremal contraction is
different from the one generally found in the literature, namely the
``contraction of an extremal ray''. However, by Proposition
\ref{P:aboutD} below, (P1) and (P2) imply that the birational map
$\f$ is the contraction of an extremal ray $R$ with respect to
$K_X+\epsilon D$, for any $0 < \epsilon$, where $R:=\mathbb R_{+}[C]$
for any curve $C \subset D$. 
\end{remark}

\begin{proposition}\label{P:aboutD}  
Let $X$ be a nonsingular projective variety of dimension $n$, and let $\varphi \colon X 
\longrightarrow~Y$ be a birational extremal contraction. Let $D$ denote 
the exceptional locus of $\f$ (with the reduced structure). Suppose 
that $\codim_X D=1$ and $\dim \f(D)=0$. Then $D$ is irreducible and 
$-D$ is  $\f$-ample. Furthermore, $\f$ is the contraction of an 
extremal ray of the cone $\overline{NE}(X)$.
\end{proposition}  
\begin{proof}  
Suppose that there are two distinct irreducible components  
$D_1$ and $D_2$ of $D$, with $\codim_X D_1=1$. Let $H_1, H_2, \dots H_{n-2} \subset X$ be general hyperplane 
sections of $X$. Let $H= \cap_{i=1}^{n-2}H_i$. Then $H$ is a smooth surface and $H \cap D_1$ is an 
irreducible curve $C_1$ on $D_1$. Then 
\[ 
(D_1 \cdot C_1)_X=(D_1|_H \cdot C_1)_H=(C_1^2)_H <0 
\] 
by the negativity of self-intersection of contractible curves. 
Now, if $C_2 \subset D_2$ is a curve 
that is not contained in $D_1$, we have that $(D_1 \cdot C_2) 
\geq 0$. But this contradicts the fact that $\f$ is extremal, because $C_1$ 
and $C_2$ cannot be numerically proportional.  Therefore $D$ is irreducible and for any curve $C \subset D$ we have 
$(D \cdot C)<0$. 
 
In order to show that $-D$ is $\f$-ample,  
we only need to show that the divisor $-D|_D$ is ample on $D$
(\cite[Theorem III.4.7.1.]{Groth3}). Then, by Kleiman's criterion of ampleness, it is enough to show that   
$(-D \cdot Z) >0$ for any $Z \in \overline{NE}(D)$.

The curves of $D$ are numerically proportional on $X$, therefore $Z 
\equiv r C_1$ on $X$ for some positive rational number $r$, since $Z\neq 0$ in $\overline{NE}(D)$. Therefore 
$(-D \cdot Z)=-r(D \cdot C_1) >0,$ 
and hence $-D$ is relatively ample and 
the divisor $-D|_D$ is ample on $D$. 
 
In order to show that $\f$ is a contraction of an extremal ray, let 
$A$ be an ample divisor on $Y$. Then, by the contraction theorem of 
extremal rays (\cite{Kenji}), the face $(\f^*A)^\perp$ of 1-cycles intersecting 
trivially with $\f^*A$ in $\overline{NE}(X)$ contains an extremal ray 
$R$ (i.e. an edge of the cone $\overline{NE}(X)$) and we have a 
contraction $cont_R \colon X \lright Y'$ of the extremal ray $R$. 
 
But any curve $C$ in $R$ is contracted by $\f$, because 
$(\f^*A \cdot C)=0$. Therefore $C$ is numerically proportional to the
curve $C_1$ in $D$, and hence the extremal ray $R$ is generated by 
$C_1$. Therefore $\f$ and $cont_R$ contract the same curves. This
implies $\f =cont_R $. 
\end{proof}  
 
Now we consider again our situation: let  $X$ be a smooth projective 
3-fold and $\f \colon X \lright Y$ a $K$-trivial birational extremal 
contraction, contracting a divisor $D$ to a point $q \in Y$.   

\subsection{Description of the exceptional divisor $D$}\label{DelP}    
By the adjunction formula,  
\begin{equation}\label{E:adjunction} 
\OO_D(-K_D) \simeq \OO_X(-(K_X +D))|_D \equiv \OO_X(-D)|_D, 
\end{equation} 
and hence $\OO_D(-K_D)$ is ample 
on $D$, because $-D$ is $\varphi$-ample. Also, $D$ has 
only Gorenstein singularities, being a  (Cartier) divisor on a 
smooth 3-fold. Therefore $D$ is a so-called 
{\em del Pezzo surface} (i.e. Gorenstein with ample 
anticanonical bundle) of degree $d=(K_D ^2)$. Note that we allow the 
del Pezzo surface $D$ to be singular. 
 
By \cite[Theorem 5.2]{Gr}, the possibilities for the exceptional 
 divisor $D$ are further restricted by its degree and singularities; $D$ is either   
\begin{enumerate}  
\item a normal and rational del Pezzo surface of degree $5 \leq d \leq 9$ 
or 
\item a non-normal del Pezzo surface of degree $d=7$ or   
\item a normal del Pezzo surface of degree $d \leq 4$ (rational for $d=4$). 
\end{enumerate} 
 
In order to obtain information about the 
normal bundle $\N_{D/X}$ of $D$ in $X$, note that the equivalences \eqref{E:adjunction} above 
also show that $\N_{D/X}$ is numerically equivalent to 
$\OO_D(K_D)$. In fact, we have
that  $\N_{D/X} \simeq
\OO_D(K_D)$ (using that the Euler 
characteristic is a numerical invariant \cite[Cor. 09]{Groth1} and
$\chi(\OO_D)=1$ for del Pezzo surfaces\cite{Wat}).

\subsection{The singularity at $q \in Y$}\label{singularity} 
 
Because $\f \colon X \longrightarrow~Y$ is  $K$-trivial and extremal, the singularity $q \in Y$ is a rational Gorenstein singularity (i.e. it is 
 Gorenstein, and $\f_* K_X =K_Y$).  According to \cite{Reid1},  to an isolated rational Gorenstein 3-fold point $q \in Y$ one can attach a 
natural number $k \geq 0$ such that  
\begin{enumerate} 
\item if $k=0$ then $q \in Y$ is a cDV point. 
\item if $k=1$ then $q \in Y$ is a hypersurface singularity that is locally of 
the form $x^2+y^3+f(y,z,t)=0$, where $f=yf_1(z,t)+f_2(z,t)$ and $f_1$ 
(respectively $f_2$) is a  sum of of monomials $z^at^b$ of degree $a+b 
\geq 4$ (respectively $\geq 6$).  
\item if $k=2$ then $q \in Y$ is a hypersurface singularity that is locally of 
the form $x^2+f(y,z,t)=0$, where $f$ is a  sum of of monomials of degree $\geq 4$. 
\item if $k \geq 3$ then $\mult_q Y=k$ and the embedding dimension of 
$q \in Y$ is $k+1$. In particular, for $k=3$, $q \in Y$ is still a hypersurface 
singularity, and for $k=4$ it is a complete intersection. 
\end{enumerate} 
\begin{remark}
\cite[Prop. (2.13)]{Reid1} implies that if the exceptional 
locus of $\f \colon X \lright Y$ is a del Pezzo surface of degree $d$, 
then the invariant $k$  is equal to $d$.
We also have that $X$ is the (weighted) blow-up of $Y$ at 
$q$ (\cite[Theorem (2.11)]{Reid1}). In particular, when $d \geq 3$, $X$ is 
the  blow-up of $Y$ at $q$. 
\end{remark}

 \section{Normal rational del Pezzo surfaces of degree $\geq 5$}\label{gen} 
\subsection{General information about del Pezzo surfaces}
 
\begin{definition} 
A 2-dimensional (possibly singular) projective variety $D$ is called a 
\underline{del Pezzo surface} if it has only Gorenstein singularities, and its 
anticanonical sheaf $\omega_D^{-1}$ is ample. We call the intersection number 
$d=(\om^2)$ the \underline{degree} of the del Pezzo surface $D$. 
\end{definition} 
Normal rational del Pezzo surfaces were classified by Hidaka and Watanabe in 
\cite{Wat}; non-normal ones by Reid in \cite{Reid}.   
In this section we enumerate some facts that will be used 
subsequently in our proof of analytic rigidity. 

Let $D$ be a normal del Pezzo surface of degree $d$ and $\pi \colon \DD \lright D$ a
minimal resolution of $D$.  Then $\DD$ is either a cone over an
elliptic curve, or $\DD$ is rational. In the latter case, $D$ is either
$\PP^2$ ($d=9$), $\PP^1 \times \PP^1$ ($d=8$), a 
singular quadric in $\PP^3$ ($d=8$), or its minimal resolution $\DD$ 
is the blow-up of $9-d$ points in almost general position on 
$\PP^2$ (\cite{Wat}). Note that if $d \geq 3$, then $\omega_D^{-1}$ is very ample and its global 
sections yield an embedding of $D$ into $\PP^d$ as a subvariety of 
degree $d$. This embedding defines a projectively normal variety
and is defined by quadric equations except for the case $d=3$
(\cite[Theorem 4.4]{Wat}). 

\begin{proposition}\cite[Prop. 4.2]{Wat}\label{P:cohom}
Let $D$ be a normal del Pezzo surface. Then
\begin{enumerate}
\item The anticanonical system $|-K_D|$ of $D$ contains a nonsingular elliptic curve.
\item $H^1(D, \OO_D(\nu K_D))=0$ for all $\nu \in \mathbb Z$.
\item If $\deg $D$=d$, then
\begin{equation*}
\dim H^0(D,\OO_D(-\nu K_D))=
\begin{cases}
\frac{(\nu +1)\nu}{2}d+1 &\text{if } \nu \geq 0\\
0 &\text{if } \nu < 0.
\end{cases}
\end{equation*}
\end{enumerate}
\end{proposition}

\subsection{Singularities of normal rational del Pezzo surfaces of
degree $\geq 5$}\label{delp} 
In the following, we present a complete description of the 
configurations of points the blow-up of which gives the minimal 
resolution of a normal rational del Pezzo surface of degree $d \geq
5$.  

\noindent{\bf Notation.} For convenience, we introduce the following
notations/definitions (\cite{Demazure}):  Let $\Sigma=\{p_1, \dots, p_r\}$ be a finite set of points on 
$\PP^2$ (infinitely near points allowed) and assume that $|\Sigma| 
\leq 8$. Denote by $\Sigma_i$ the subsystem $\{p_1, \dots, p_i\}$, $1 
\leq i \leq r$, and let $V(\Sigma_i) \lright \PP^2$ be the blow-up 
of $\PP^2$ with center $\Sigma_i$. Then there exists a sequence of 
blow-ups 
\[ 
V(\Sigma)=V(\Sigma_n)\lright V(\Sigma_{n-1})\lright \dots \lright 
V(\Sigma_1)\lright \PP^2. 
\] 
Denote by $E_i \subset V(\Sigma_i)$ the exceptional locus of the morphism 
$V(\Sigma_i)\lright V(\Sigma_{i-1})$. 
\begin{definition}\cite[III.2]{Demazure} 
The points of $\Sigma$ are in \underline{general} (respectively \underline{almost general}) \underline{position} 
if 
\begin{enumerate} 
\item No three (resp. four) of them are collinear. 
\item No six (resp. seven) of them are on a conic. 
\item All the points are distinct (resp. for all $1 \leq i \leq r$, 
the point $p_{i+1} \in V(\Sigma_i)$ does not lie on any proper 
transform $\hat E_j$ of some $E_j$, $1 \leq j \leq i$, with $(\hat 
E_j^2)=-2$).
\item When $|\Sigma| = 8$, there exists no singular cubic which passes 
through all the points of $\Sigma$ and has one of them as the singular 
point (no corresponding condition in the almost general case). 
\end{enumerate} 
\end{definition} 
 
Let $D$ be a normal rational del Pezzo surface of degree $d \geq 5$
and let $\pi \colon \DD \lright D$ be its minimal resolution. Suppose
$D$ is not $\PP^1 \times \PP^1$ or a singular quadric in $\PP^3$. Then
$\DD \simeq V(\Sigma)$ for some set $\Sigma$ of $9-d \leq 4$ points
in almost general position on $\PP^2$.

Up to a projective automorphism of $\PP^2$ (and its extensions to the
blow-up spaces), there are 22 different configurations of at most 4 
points in almost general position (including $\Sigma = \emptyset$,
when $D \simeq \PP^2$). Below we describe these configurations (giving
a representative for each configuration) and the corresponding blow-ups.
 
Fix the homogeneous coordinates $[x_0 \colon x_1 \colon x_2]$ on 
$\PP^2$. Then on the affine open $U_0=\{x_0 \neq 0\} \simeq \mathbb A^2$ we have local 
coordinates $x: = \frac{x_1}{x_0}$ and $y := \frac{x_2}{x_0}$.  In the following, all blow-ups will be of 
(possibly infinitely near) points of $U_0 \subset \PP^2$. For any
point $p \in U_0$, we have the natural coordinates on the
exceptional locus $\EE$ of the blow-up $\Bl_p U_0$, using $\EE
\simeq \PP(T_p U_0) \simeq \PP(\mathbb A^2) \simeq \PP^1$.  For a point $q \in \Bl_p U_0$ we write $q
\lright p$ if $q$ is infinitely near to $p$, i.e. if $q \in \EE$.

\subsubsection{Blowing up one point} 
 
\chunkk Let $\Sigma_1=\{p_1\}$, where $p_1 \in \PP^2$ and let $\s_1 \colon 
S_1=V(\Sigma_1) \lright \PP^2$ be the blow-up of $p_1$ on $\PP^2$. Up
to a projective automorphism of $\PP^2$, we can assume $p_1=[1:0:0]$.
 
\subsubsection{Blowing up  two points} 
 
\chunkk Let $\Sigma_2=\Sigma_1 \cup \{p_2\}$, where $p_2\in \PP^2$ and $p_2 \neq p_1$. Let $\s_2 \colon 
S_2=V(\Sigma_2) \lright S_1$ be the blow-up of $p_2$ on $S_1$. Up
to a projective automorphism of $\PP^2$, we can assume $p_1=[1:0:0]$
and $p_2=[1:0:1]$.
 
\chunkk Let $\Sigma_7=\Sigma_1 \cup \{p_2\}$, where $p_2 \lright p_1$. Let $\s_7 \colon 
S_7=V(\Sigma_7) \lright S_1$ be the blow-up of $p_2$ on $S_1$. Up
to a projective automorphism of $\PP^2$ (and its lifting to a
projective automorphism of $V(\Sigma_1)$), we can assume $p_1=[1:0:0]$
and $p_2=[1:0] \in \PP(T_{p_1} \PP^2 ) \subset V(\Sigma_1)$.
 
\subsubsection{Blowing up three points} 
 
\case{The three points are on $\PP^2$} 
 
\chunkk Let $\Sigma_3=\Sigma_2 \cup \{p_3\}$, where $p_3 \in 
\PP^2$ and $p_1,p_2,p_3$ are not collinear. Let $\s_3 \colon 
S_3=V(\Sigma_3) \lright S_2$ be the blow-up of $p_3$ on $S_2$. We can
assume $p_1=[1:0:0]$, $p_2=[1:0:1]$ and $p_3=[1:1:0]$ .

\chunkk Let $\Sigma_4=\Sigma_2 \cup \{p_3\}$, where $p_3 \in \PP^2$ and $p_1, p_2, p_3$ are collinear. Let $\s_4 \colon 
S_4=V(\Sigma_4) \lright S_2$ be the blow-up of $p_3$ on $S_2$. We can
assume $p_1=[1:0:0]$, $p_2=[1:0:1]$ and $p_3=[1:0:-1]$ .
 
\case{Two of the three points are on $\PP^2$} 
\chunkk Let $\Sigma_5=\Sigma_2 \cup \{p_3\}$, where $p_3 \lright p_1$ and such that $p_1, p_2, p_3$ are not collinear. Let $\s_5 \colon 
S_5=V(\Sigma_5) \lright S_2$ be the blow-up of $p_3$ on $S_2$. We can
assume $p_1=[1:0:0]$, $p_2=[1:0:1]$ and $p_3=[1:0] \in \PP(T_{p_1} \PP^2) \subset V(\Sigma_2)$.
 
\chunkk Let $\Sigma_6=\Sigma_2 \cup \{p_3\}$, where $p_3 \lright p_1$, such that $p_1, p_2, p_3$ are collinear. Let $\s_6 \colon 
S_6=V(\Sigma_6) \lright S_2$ be the blow-up of $p_3$ on $S_2$. We can
assume $p_1=[1:0:0]$, $p_2=[1:0:1]$ and $p_3=[0:1] \in \PP(T_{p_1}
\PP^2 ) \subset V(\Sigma_2)$.

\case{Only one of the three points is on $\PP^2$}  
\chunkk Let $\Sigma_8=\Sigma_7 \cup \{p_3\}$, where $p_3  \lright p_2
\lright p_1$,  and $p_1, 
p_2, p_3$ are not collinear. Let $\s_8 \colon 
S_8=V(\Sigma_8) \lright S_7$ be the blow-up of $p_3$ on $S_7$. We can assume $p_1=[1:0:0]$, $p_2=[1:0] \in \PP(T_{p_1} \PP^2) \subset V(\Sigma_1)$ and
$p_3=[1:1] \in \PP(T_{p_2}  V(\Sigma_1)) \subset V(\Sigma_7)$.

\chunkk Let $\Sigma_9=\Sigma_7 \cup \{p_3\}$, where $p_3  \lright p_2
\lright p_1$ and $p_1, 
p_2, p_3$ are collinear. Let $\s_9 \colon 
S_9=V(\Sigma_9) \lright S_7$ be the blow-up of $p_3$ on $S_7$. We can
assume $p_1=[1:0:0]$,  $p_2=[1:0] \in \PP(T_{p_1} \PP^2) \subset V(\Sigma_1)$ and
$p_3=[1:0] \in \PP(T_{p_2}  V(\Sigma_1)) \subset V(\Sigma_7)$.
 
\subsubsection{Blowing up four points} 
 
\case{All four points are on $\PP^2$}
\chunkk Let $\Sigma_3'=\Sigma_3 \cup \{p_4\}$, where $p_4  \in \PP^2$, with no
three points collinear. Let $\s_3' \colon 
S_3'=V(\Sigma_3') \lright S_3$ be the blow-up of $p_4$ on $S_3$. We can
assume $p_1=[1:0:0]$, $p_2=[1:0:1]$, $p_3=[1:1:0]$  and $p_4=[1:1:-1]$ .

\chunkk Let $\Sigma_4'=\Sigma_4 \cup \{p_4\}$, where $p_4 \in 
\PP^2$, with
$p_1$, $p_2$ and $p_3$ collinear. Let $\s_4' \colon 
S_4'=V(\Sigma_4') \lright S_4$ be the blow-up of $p_4$ on $S_4$. We can
assume $p_1=[1:0:0]$, $p_2=[1:0:1]$, $p_3=[1:0:-1]$  and $p_4=[1:1:0]$ .

 \case{Three of the points are on $\PP^2$}
 
\chunkk Let $\Sigma_3''=\Sigma_3 \cup \{p_4\}$, where $p_4  \lright p_1$ and such that no
three points are collinear. Let $\s_3'' \colon 
S_3''=V(\Sigma_3'') \lright S_3$ be the blow-up of $p_4$ on $S_3$. We can
assume $p_1=[1:0:0]$, $p_2=[1:0:1]$, $p_3=[1:1:0]$  and $p_4=[1:0] \in \PP(T_{p_1} \PP^2 ) \subset V(\Sigma_3)$.

\chunkk Let $\Sigma_3'''=\Sigma_3 \cup \{p_4\}$, where $p_4  \lright p_1$ and $p_1, p_2, p_4$ are collinear. Let $\s_3''' \colon 
S_3'''=V(\Sigma_3''') \lright S_3$ be the blow-up of $p_4$ on $S_3$. We can
assume $p_1=[1:0:0]$, $p_2=[1:0:1]$, $p_3=[1:1:0]$  and $p_4=[0:1] \in \PP(T_{p_1} \PP^2 ) \subset V(\Sigma_3)$.

\chunkk Let $\Sigma_4''=\Sigma_4 \cup \{p_4\}$, where $p_4 \lright
p_1$ and is not collinear with $p_1,p_2,p_3$. Let $\s_4'' \colon 
S_4''=V(\Sigma_4'') \lright S_4$ be the blow-up of $p_4$ on $S_4$. We can
assume $p_1=[1:0:0]$, $p_2=[1:0:1]$, $p_3=[1:0:-1]$  and $p_4=[1:0]
\in \PP(T_{p_1} \PP^2 ) \subset V(\Sigma_4)$.

\case{Two of the points are on $\PP^2$}

\chunkk Let $\Sigma_5'=\Sigma_5 \cup \{p_4\}$, where $p_4  \lright p_2$ and such that no
three points are collinear. Let $\s_5' \colon 
S_5'=V(\Sigma_5') \lright S_5$ be the blow-up of $p_4$ on $S_5$. We can
assume $p_1=[1:0:0]$, $p_2=[1:0:1]$, $p_3=[1:0] \in \PP(T_{p_1} \PP^2)
\subset V(\Sigma_2)$ and $p_4=[1:0] \in \PP(T_{p_2} \PP^2 ) \subset V(\Sigma_5)$.

\chunkk Let $\Sigma_5''=\Sigma_5 \cup \{p_4\}$, where $p_4  \lright p_3
\lright p_1$ and such that no
three points are collinear. Let $\s_5'' \colon 
S_5''=V(\Sigma_5'') \lright S_5$ be the blow-up of $p_4$ on $S_5$. We can
assume $p_1=[1:0:0]$, $p_2=[1:0:1]$,  $p_3=[1:0] \in \PP(T_{p_1}
\PP^2) \subset V(\Sigma_2)$ and $p_4=[1:1] \in \PP(T_{p_3}
(V(\Sigma_2)) \subset V(\Sigma_5)$.

\chunkk Let $\Sigma_5'''=\Sigma_5 \cup \{p_4\}$, where $p_4  \lright p_3
\lright p_1$  and $p_1, p_3, p_4$ are collinear. Let $\s_5''' \colon 
S_5'''=V(\Sigma_5''') \lright S_5$ be the blow-up of $p_4$ on $S_5$. We can
assume $p_1=[1:0:0]$, $p_2=[1:0:1]$,  $p_3=[1:0] \in \PP(T_{p_1}
\PP^2) \subset V(\Sigma_2)$ and $p_4=[1:0] \in \PP(T_{p_3}
(V(\Sigma_2)) \subset V(\Sigma_5)$.

\chunkk Let $\Sigma_6'=\Sigma_6 \cup \{p_4\}$, where $p_4  \lright p_2$. Let $\s_6' \colon 
S_6'=V(\Sigma_6') \lright S_6$ be the blow-up of $p_4$ on $S_6$. We can
assume $p_1=[1:0:0]$, $p_2=[1:0:1]$,  $p_3=[0:1] \in \PP(T_{p_1}
\PP^2 ) \subset V(\Sigma_2)$ and $p_4=[1:0] \in \PP(T_{p_2} \PP^2 ) \subset V(\Sigma_6)$.
 
\chunkk Let $\Sigma_6''=\Sigma_6 \cup \{p_4\}$, where $p_4  \lright p_3
\lright p_1$ and $p_1, p_3, p_4$ are not collinear. Let $\s_6'' \colon 
S_6''=V(\Sigma_6'') \lright S_6$ be the blow-up of $p_4$ on $S_6$. We can
assume $p_1=[1:0:0]$, $p_2=[1:0:1]$,  $p_3=[0:1] \in \PP(T_{p_1}
\PP^2 ) \subset V(\Sigma_2)$ and $p_4=[1:1] \in \PP(T_{p_3}
(V(\Sigma_2)) \subset V(\Sigma_6)$.
 
\case{One point is on $\PP^2$}
 
\chunkk Let $\Sigma_8'=\Sigma_8 \cup \{p_4\}$, where $p_4  \lright
p_3 \lright p_2 \lright p_1$, and such that no
three points are collinear. Let $\s_8' \colon 
S_8'=V(\Sigma_8') \lright S_8$ be the blow-up of $p_4$ on $S_8$. We can assume $p_1=[1:0:0]$, $p_2=[1:0] \in \PP(T_{p_1} \PP^2) \subset V(\Sigma_1)$,
$p_3=[1:1] \in \PP(T_{p_2}  V(\Sigma_1)) \subset V(\Sigma_7)$ and
$p_4=[1:0] \in \PP(T_{p_3} (V(\Sigma_7)) \subset V(\Sigma_8)$.

\chunkk Let $\Sigma_9'=\Sigma_9 \cup \{p_4\}$, where $p_4  \lright
p_3 \lright p_2 \lright p_1$, such that  $p_1, 
p_2, p_3$ are collinear. Let $\s_9' \colon 
S_9'=V(\Sigma_9') \lright S_9$ be the blow-up of $p_4$ on $S_9$. We
can assume $p_1=[1:0:0]$, $p_2=[1:0] \in \PP(T_{p_1} \PP^2) \subset
V(\Sigma_1)$, 
$p_3=[1:1] \in \PP(T_{p_2}  V(\Sigma_1)) \subset V(\Sigma_7)$ and
$p_4=[1:0] \in \PP(T_{p_3} (V(\Sigma_7)) \subset V(\Sigma_9)$.
 
\begin{remark}\label{R:config}
Observe that, given a set of points $\Sigma$ in almost general 
position on $\PP^2$, with $|\Sigma| \leq 4$, there exists a projective automorphism of $\PP^2$ 
that transforms $\Sigma$ into one of the point-sets listed above ($\Sigma_i$, $\Sigma_i'$, $\Sigma_i''$), say $\widetilde 
\Sigma$. We say that $\Sigma$ and $\widetilde 
\Sigma$ have the same configuration, and we have $V(\Sigma) \simeq V(\widetilde\Sigma)$. 
\end{remark}

\begin{corollary}\label{C:sing} 
A normal rational del Pezzo surface $D$ of degree $d \geq 5$ is either
$\PP^2$, $\PP^1 \times \PP^1$ or a 
singular quadric in $\PP^3$, or it can have the 
following singularities: 
\begin{enumerate} 
\item No singularities if $\DD$ is $S_1$, $S_2$, $S_3$ or $S_3'$. 
\item One $A_1$-singularity in the case $\DD$ is $S_4$,  $S_5$, $S_7$,
$S_3''$ or $S_4'$. 
\item Two $A_1$-singularities in the case $\DD$ is  $S_6$, $S_5'$ or 
$S_3'''$. 
\item One $A_2$-singularity in the case $\DD$ is  $S_4''$, $S_8$ or $S_5''$. 
\item One $A_1$- and one $A_2$-singularity in the case $\DD$ is $S_6'$, 
$S_9$ or $S_5'''$. 
\item One $A_3$-singularity in the case $\DD$ is  $S_6''$ or $S_8'$. 
\item One $A_4$-singularity in the case $\DD$ is  $S_9'$. 
\end{enumerate} 
\end{corollary}

\section{Analytic rigidity, formal equivalence and infinitesimal extensions}\label{first} 
 
\subsection{Analytic rigidity and formal equivalence}
A standard tool for showing analytic equivalence is a criterion due 
to Grauert \cite{Grauert} and Hironaka-Rossi 
\cite{Hir} that reduces the problem of showing analytic equivalence of
embeddings to 
that of showing formal 
equivalence (\cite[Theorem 3]{Hir}). We also use the following lemma (see
\cite[Lemma 9]{Hir}):

\begin{lemma} \label{L:Hiron}  
Let $X$ be a nonsingular complex manifold of dimension $n$ and let $D$ be a reduced complex subspace with ideal sheaf $\mathcal I$ (i.e. $\sqrt {\mathcal I}={\mathcal I}$). Suppose that there exists an integer $\nu _0 \geq 2$ such that  
 $H^1(D, \mathcal T_X \otimes \mathcal I ^{\nu}/\mathcal I ^{\nu +1})=0$ for  
any $\nu \geq \nu _0$. Then a $\nu$-equivalence $(\nu \geq \nu _0)$ of 
 $D$ with a complex subspace $D'$ of 
a complex space $X'$, where $X'$ has the same dimension as $X$ at any point 
 of $D'$, extends to a formal equivalence.  
\end{lemma}  
  
By \underline{$\nu$-equivalence} we mean an isomorphism of complex
 spaces $(D, \mathcal O_X/\mathcal I_D ^{\nu}) \simeq (D', \mathcal
 O_{X'}/\mathcal I_{D'} ^{\nu})$, where $\I_D$, respectively $\I_{D'}$
 is the ideal sheaf of $D$ in $X$, respectively of $D'$ in
 $X'$. \underline{Formal equivalence} means an isomorphism of the
 completions:
 $$\hat X=\varprojlim (D, \mathcal O_X/\mathcal I_D ^{\nu})
\simeq \hat{X'} =\varprojlim (D', \mathcal O_{X'}/\mathcal {I_{D'}}
^{\nu}).$$

\begin{remark}   We extracted the condition
$H^1(D, \mathcal T_X \otimes \mathcal I ^{\nu}/\mathcal I ^{\nu
+1})=0$, $\forall \nu \geq \nu_0$  from the proof of \cite[Lemma 9]{Hir} (for any $\nu$, the obstruction to extending a
$\nu$-equivalence to a $(\nu +1)$-equivalence lies in the above cohomology group). In the original version $D$ was assumed to 
be the exceptional locus of the blowing up of a point of a complex space.  
 While in our case the del Pezzo surface $D$ is the exceptional locus 
 of the blowing up of $q$ on $Y$, we want to use this criterion with 
 $\nu_0 =2$, and hence we need to have a better control 
 on $\nu_0$ than the proof of \cite[Lemma 9]{Hir} provides. The cohomological condition of the vanishing of 
 $H^1(D, \mathcal T_X \otimes \mathcal I ^{\nu}/\mathcal I ^{\nu +1})$ gives us this control. 
\end{remark}   

In order to use Lemma \ref{L:Hiron}, we need to establish a criterion of when 
two different embeddings of a 2-dimensional variety $D$ into smooth three-folds are 
$2$-equivalent. In the following we 
give a criterion for 2-equivalence using the theory of infinitesimal extensions of $\OO_D$. While there is a general theory of infinitesimal extensions of
schemes, there seems to be no appropriate reference for our purpose.  Therefore
in \S \ref{2str} we give a self-contained presentation in the
case of infinitesimal extensions of a 2-dimensional variety with
isolated hypersurface singularities that can be embedded into a smooth
3-fold.

\subsection{2-structures on $D$}\label{2str}

Let $D$ be a projective 2-dimensional variety over an algebraically 
closed field $k$. 
(We will later use the results of this section with $k=\C$ in 
establishing analytic rigidity.) 
Suppose $D$ has only isolated 
hypersurface-singularities. We assume $D$ can be embedded into a smooth projective 3-fold.    
 
Fix such an embedding $D \subset X$ and let $\I_D$ denote the ideal 
sheaf of $D$ in $X$. The conormal bundle of $D$ in $X$ is 
$\I_D/\I_D^2$. Then we have the following exact sequence  
\[ 
0 \lright \II \lright \OO_X/\I_D^2 \lright \OO_D \lright 0 .
\] 
Therefore $2D=(D,\OO_X/\I_D^2)$ is an {\it infinitesimal 
extension} of $D$ by the sheaf $\LL=\I_D/\I_D^2$ 
(i.e. $\LL$ can be considered as an ideal sheaf with square 0 on the 
scheme $2D$, with $\OO_{2D}/\LL \simeq \OO_D$ \cite[Ex. II.8.7]{Hartsh}).

In the following sections we describe the 
infinitesimal extensions of 
$D$ by $\II$, and 
give a criterion to establish when two embeddings of $D$ into smooth projective 3-folds, with isomorphic 
normal bundles, are 2-equivalent. 
 
\subsubsection{Local computations}\label{ext-local} 
 
Let $U=\Spec A$ be an affine open in X, and suppose $D \cap U$ is 
given by the ideal $I$. Denote $B \: = A/I$. Then we have an exact 
sequence 
\[ 
0 \lright I/I^2 \lright A/I^2 \lright B \lright 0 
\] 
which gives an infinitesimal extension of the $k$-algebra $B$ by the 
$B$-module $I/I^2$.

Let us review some facts about infinitesimal extensions 
of $k$-algebras, based on \cite[Sect. 18]{Groth} (here $k$ is an 
arbitrary field). 
 
\begin{definition}\cite{Groth} 
Let $B$ be a $k$-algebra. An \underline{infinitesimal extension of 
the $k$-algebra $B$ by a $B$-module $L$}  
is a $k$-algebra $E$ such that  
$L$ can be embedded in 
$E$ as an ideal whose square is 0 and such that $E/L \simeq B$. 
\end{definition} 
 
\begin{definition} 
Two infinitesimal extensions $E$ and $E'$ of $B$ by $L$ are said to be \underline{isomorphic} if there is a $k$-algebra homomorphism $E \lright E'$ which makes the following diagram commutative: 
\[ 
\xymatrix{0 \ar[r] & L \ar@{=}[d] \ar[r] & E \ar[d] \ar[r] & B \ar@{=}[d] \ar[r] & 0\\ 
0 \ar[r] & L \ar[r] & E' \ar[r] & B \ar[r] & 0 } 
\] 
\end{definition}  
In the spirit of \cite{Groth}, we have the following construction: 
\begin{definition} Let $0  \lright L\xrightarrow{\al} E
\xrightarrow{\bb} B \lright 0$ be an infinitesimal extension of $B$ by $L$ and let $L'$ be 
another $B$-module. Suppose we have a $B$-module homomorphism $\f 
\colon L \lright L'$. We define the \underline{push-out} $E \oplus _{\f} L'$ of 
$E$ and $L'$ under $\f$ by $E \oplus _{\f}~L'= \coker ((\bb,-\f) 
\colon L \lright E \oplus L')=(E \oplus L')/\{(\bb(l), -\f (l)) : l 
\in~L\}$. The ring structure on $E \oplus _{\f}~L'$ is given by $(e,l') \cdot 
(e',l'')=(ee',el''+e'l')$. Then $0 \lright L \xrightarrow{\al_\f} E \oplus _{\f} L'  \xrightarrow{\bb_\f} B \lright 0$ is an infinitesimal extension of $B$ by 
$L'$ and we have the following commutative diagram: 
\[ 
\xymatrix{0 \ar[r] & L \ar[d]_{\f} \ar[r]^\al & E \ar[d]_{g_\f} \ar[r]^\bb & B \ar@{=}[d] \ar[r] & 0\\ 
0 \ar[r] & L' \ar[r]_{\al_\f \quad} & E \oplus _{\f} L' \ar[r]_{\quad \bb_\f} & B 
\ar[r] & 0 } 
\] 
where $\al_\f (f)=(0,f)$, $\bb_\f(a,f)=\bb (a)$ and $g_\f (a)=(a,0)$. 
\end{definition}

Let $\Exal_k (B,L)$ denote the set of isomorphism classes of 
infinitesimal extensions of $B$ by $L$ . Then it has a 
$k$-module structure (\cite{Groth}); in particular, 
given an infinitesimal extension $\EE$ and a constant $\lambda 
\in k$, the multiplication $\lambda \EE$ is defined by the 
infinitesimal extension $ 0 \lright L \lright E \oplus _{\lambda} L 
\lright B \lright 0$, where $E \oplus _{\lambda} L= (E \oplus 
L)/\{(l,-\lambda l) : l \in L\}$. Note that multiplication by $0 \in k$ gives the trivial infinitesimal extension $B \oplus L$.

\begin{lemma} 
Let $D$ be a projective 2-dimensional variety over an algebraically 
closed field $k$. Suppose $D$ can be embedded into a smooth projective 
3-fold $X$. Let $U=\Spec A$ be an affine open in X, and suppose $D 
\cap U=\Spec B$, where $B=A/I$. Then for any $B$-module $L$, there is an isomorphism $$\Exal_k(B,L) \simeq 
\Ext_B^1(\Omega_{B/k}, L).$$ 
\end{lemma} 
\begin{proof} 
By \cite[Theorem 4.2.2]{L-S}, $\Exal_k(B,L) \simeq T^1(B/k, L)$, where $T^1$ is the first 
cotangent functor. Also, \cite[Lemma 3.1.1]{L-S} implies that if $P$ 
is a polynomial ring over $k$ with $B \simeq P/J$ ($J$ an ideal of 
$P$), then we have the isomorphism 
\[ 
T^1(B/k, L)=\coker(\Hom_P(\Omega_{P/k},L) \lright \Hom(J/J^2,L)). 
\] 
 
However, $D$ is reduced and it is a Cartier divisor in the nonsingular 
3-fold $X$, and hence $\Spec B= D \cap U$ is a reduced local complete 
intersection in $U=\Spec A$. The notion of $\Spec B$ being a local complete 
intersection is independent of the nonsingular variety containing 
it, hence $\Spec B$ is also a reduced local complete 
intersection in $\Spec P$. Therefore, by \cite[Ex. 16.17]{Eisenb} we obtain that the sequence 
\begin{equation}
0 \lright J/J^2 \lright \Omega_{P/k} \otimes_P B \lright \Omega_{B/k} \lright 0 
\end{equation} 
is exact, and so
$\coker(\Hom_P(\Omega_{P/k},L) \lright \Hom(J/J^2,L)) \simeq 
\Ext_B^1(\Omega_{B/k},L). $
\end{proof} 
In our case (of infinitesimal extensions of $\OO_D$ by $\II$), $L=I/I^2$.  In the following discussions we use that $\Exal_k(B,I/I^2) \simeq \Ext^1_B(\Omega_{B/k}, I/I^2)$. 
\begin{remark} If the open set $U$ does not contain any singular points of $D$, then 
$\Ext^1_B(\Omega_{B/k}, I/I^2) =0$, and therefore  $\text 
{Exalcom}_k(B,I/I^2)$ only contains the trivial infinitesimal extension $B \oplus I/I^2 \simeq A/I^2$. 
\end{remark}

\subsubsection{General discussion about infinitesimal extensions of $B$ by $I/I^2$}\label{ext:general}

We keep the notations from \S \ref{ext-local}. As $D$ is a 
local complete intersection in the smooth 3-fold $X$, we have the 
exact sequence 
\begin{equation}\label{Eq:exact} 
0 \lright I/I^2 \lright \Omega_{A/k} \otimes_A B \lright \Omega_{B/k} \lright 0. 
\end{equation} 
From here, we obtain the long exact sequence 
\[ 
\begin{aligned} 
0 \lright &\Hom_B(\Omega_{B/k},I/I^2) \lright \Hom_B(\Omega_{A/k} \otimes_A B,I/I^2) \lright \Hom_B(I/I^2,I/I^2) \lright \\ 
&\Ext^1_B(\Omega_{B/k},I/I^2) \lright \Ext^1_B(\Omega_{A/k} \otimes_A B,I/I^2)=0, 
\end{aligned} 
\] 
where the last term vanishes because $\Omega_{A/k}$ is locally 
free. Therefore we obtain that ${\Exal _k}(B,I/I^2) \simeq \Ext^1_B(\Omega_{B/k},I/I^2)$ is a homomorphic image of $\Hom_B(I/I^2,I/I^2)$. The correspondence goes as follows: let $\f \in \Hom_B(I/I^2,I/I^2)$ and consider the standard exact sequence 
\[ 
\xymatrix{0 \ar[r] & I/I^2 \ar[r]^\al & A/I^2 \ar[r]^\bb & B \ar[r] & 0}. 
\] 
 
Then the infinitesimal extension corresponding to $\f$ is given by $$\A_\f : 0 \lright I/I^2 \lright A_\f \lright B \lright 0,$$ where $A_\f=A/I^2 \oplus _{\f} I/I^2 $. We obtain the following commutative diagram: 
\begin{equation}\label{CD:A_fi} 
\xymatrix{\A_0 : & 0 \ar[r] & I/I^2 \ar[d]_{\f} \ar[r]^\al & A/I^2 \ar[d]_{g_\f} \ar[r]^\bb & B \ar@{=}[d] \ar[r] & 0\\ 
\A_\f: & 0 \ar[r] & I/I^2 \ar[r]_{\al_\f} & A_\f \ar[r]_{\bb_\f} & B \ar[r] & 0 } 
\end{equation} 
where $\al_\f (f)=(0,f)$, $\bb_\f(a,f)=\bb (a)$ and $g_\f (a)=(a,0)$.  
 
\begin{lemma} \label{L:deriv}
Two extensions $\A_\f$ and $\A_{\f '}$ are isomorphic if 
and only if $\f$ and $\f '$ differ by a $k$-derivation of $A/I^2$ into 
$I/I^2$.  
\end{lemma} 
\begin{proof}If $\psi \in \Der_k(A/I^2,I/I^2) \simeq \Hom_B(\Omega_{A/k} \otimes_A B,I/I^2)$ and $\f '=\f +\psi \circ \al$, then we have a morphism (and hence isomorphism) of extensions 
\[ 
\xymatrix{0 \ar[r] & I/I^2 \ar@{=}[d] \ar[r]^{\al_\f} & A_\f \ar[d]_{\Psi} \ar[r]^{\bb_\f} & B \ar@{=}[d] \ar[r] & 0\\ 
0 \ar[r] & I/I^2 \ar[r]_{\al_{\f '}} & A_{\f '} \ar[r]_{\bb_{\f '}} & B \ar[r] & 0 } 
\] 
where $\Psi (a,f)=(a,f+\psi(a))$. 
 
Conversely, if we have a morphism $\Psi$ that makes the above diagram 
commutative, then 
$\Psi(a,0)-(a,0)=(0,\psi 
(a))$ because $\bb_{\f'} \circ \Psi=\bb_{\f}$ and $\Psi \circ 
\al_{\f}=\al_{\f '}$
, where 
$\psi \colon A/I^2 \lright I/I^2$. The multiplication rule on $A_{\f 
'}$ and the fact that $\Psi$ is a ring 
homomorphism implies then that  
$\psi \in \Der_k(A/I^2,I/I^2)$. It follows then that $\Psi (a,f)=(a,f+\psi(a))$.
\end{proof} 
\begin{remark} Note that $\A_0$ is the trivial (split) extension and
$\A_1$ is the ``standard'' conormal extension $0 \lright I/I^2
\xrightarrow{\al} A/I^2 \xrightarrow{\bb} B \lright 0$. Also, given
any infinitesimal extension $\EE:  0 \lright I/I^2 \xrightarrow{u} E
\xrightarrow{v} B \lright 0$ of $B$ by $I/I^2$, there exists a $\f \in
\Hom_B(I/I^2,I/I^2)$ such that the extensions $\A_\f$ and $\EE$ are
isomorphic. (This follows from the fact that ${\Exal _k}(B,I/I^2)$ is the homomorphic image of $\Hom_B(I/I^2,I/I^2)$.) 
 \end{remark}
\begin{remark} Observe that $\f$ is isomorphism if and only if $g_\f$ is, and in that case $\A_\f$ gives rise to an extension $\f ^* \A_\f: 0 \lright I/I^2 \xrightarrow{\al_\f \circ \f} A_\f \xrightarrow{\bb_\f} B \lright 0$ isomorphic to $\A_1$ via 
\[ 
\xymatrix{0 \ar[r] & I/I^2 \ar@{=}[d] \ar[r]^{\al} & A/I^2 \ar[d]^{g_\f} \ar[r]^{\bb} & B  \ar@{=}[d] \ar[r] & 0\\ 
0 \ar[r] & I/I^2  \ar[r]_{\al_\f \circ \f} & A_\f  \ar[r]_{\bb_\f} & B \ar[r] & 0 .} 
\] 
\end{remark} 
 
\begin{lemma}\label{L:trivext} 
Suppose we have two isomorphisms of extensions 
\[ 
\xymatrix{0 \ar[r] & I/I^2 \ar@{=}[d] \ar[r]^\al & A/I^2 \ar[d]_{g\,\,} \ar@<-1ex>[d]^{\,\,\, g^\prime} \ar[r]^\bb & B \ar@{=}[d] \ar[r] & 0\\ 
0 \ar[r] & I/I^2  \ar[r]^u & E  \ar[r]^v & B \ar[r] & 0.} 
\] 
Then $g$ and $g'$ differ by a derivation  of $B$ into $I/I^2$. More precisely, there exists a $\psi \in \Der _k(B,I/I^2)$ such that $g=g'+u \circ \psi \circ \bb$. 
\end{lemma} 
\begin{proof} 
Let $\Phi := g^{-1} \circ g'$. Then we have the following commutative diagram 
\[ 
\xymatrix{0 \ar[r] & I/I^2 \ar@{=}[d] \ar[r]^{\alpha} & A/I^2 \ar[d]_{\Phi} \ar[r]^{\beta} & B \ar@{=}[d] \ar[r] & 0\\ 
0 \ar[r] & I/I^2 \ar[r]^{\alpha} & A/I^2 \ar[r]^{\beta} & B \ar[r] & 0 .} 
\] 
 
From $\beta \circ \Phi =\beta$ follows that for any $a \in A/I^2$,
$\Phi (a) =a+\alpha (f_a)$ for some unique $f_a \in I/I^2$. Using that
$\Phi$ is a $k$-algebra isomorphism, we obtain that the map $\psi :
A/I^2 \lright I/I^2$ given by $\psi(a) : =f_a$ is a $k$-derivation of
$A/I^2$, and from $\Phi \circ \alpha =\alpha$ it follows that $\psi$
factors through $B=A/I$. So we may assume $\psi \in \Der _k (B, I/I^2)
\simeq \Hom_B (\Omega _B, I/I^2)$.  Therefore we have $\Phi = id +\alpha \circ \psi \circ \beta$, and hence $g' =g + u \circ \psi \circ \bb$. 
\end{proof}

\begin{proposition}\label{P:completion}   
Let $p$ be a maximal ideal of $B$ and suppose that the completion of 
$B$ at $p$ is $\hat B=k [[x,y,z]]/(f)$ where $f$ is of order $\geq 
1$. Then $f\hat B/f^2\hat B \simeq \hat B$ and $$\Exal_k(\hat B, \hat B)\simeq \Ext ^1_{\hat 
B}(\Omega_{\hat B}, \hat B)=k [[x,y,z]]/(f,\frac{\partial f}{\partial
x},\frac{\partial f}{\partial y},\frac{\partial f}{\partial z}).$$
Moreover,  for any $\f \in \Hom ( \hat B, \hat B)$,  the corresponding infinitesimal extension $\hat \A_\f$ 
of $\hat B$ by $\hat B$ is given by $$0 \lright \hat B \xrightarrow{\hat \al_\f} \hat A_\f =k [[x,y,z,T]]/(f^2, 
T^2, \f T-f) \xrightarrow{\hat \bb_\f} \hat B \lright 0,$$ where we 
consider $\f$ as an element of $\hat B$. The maps are defined by $\hat\al_\f (1)=T$ and $\hat\bb_\f(T)=0$.  
\end{proposition} 
\begin{proof} 
The only thing needing attention is the isomorphism $$\hat A_\f 
\simeq k[[x,y,z,T]]/(f^2, T^2, \f T-f).$$ Using the definition of
$\hat A_\f$ as $k [[x,y,z]]/(f^2) \oplus _\f~\hat B$, the isomorphism is given by the map $$\Phi \colon k 
[[x,y,z]]/(f^2) \oplus _\f \hat B \lright 
k[[x,y,z,T]]/(f^2, T^2, \f T-f),$$ defined by $\Phi(h,g)=h+Tg$. 
\end{proof} 
\begin{remark} 
Note that for any invertible $\f$, we have $k[[x,y,z]]/(f^2) \simeq \hat A_\f$ via $\hat g_\f$. 
\end{remark} 

We are looking at infinitesimal extensions of $D$ obtained from 
   embedding $D$ into smooth 3-folds. Therefore, locally, the 
   extensions we are concerned about need to have embedding dimension 
   3. From Proposition \ref{P:completion}, we obtain the following: 
\begin{corollary}\label{C:embdim} 
Under the conditions of Proposition \ref{P:completion}, $\hat A_\f$ has embedding dimension 3 if and only if $\f$ is invertible. 
\end{corollary} 

Now suppose that $U$ contains one (and only one) singular point $p$ 
of $D$. In this case we have the following: 
\begin{proposition}\label{P:embdim3} 
Under the assumptions of Proposition \ref{P:completion}, an infinitesimal extension $A_\f$ of $B$ by $I/I^2$ has embedding 
dimension 3 at $p$       
if and only if $\f$ is isomorphism; in that case $g_\f:A/I^2 \lright A_\f$ is an isomorphism. 
\end{proposition} 
\begin{proof} 
We look again at the definition of $A_\f$ via \eqref{CD:A_fi}: 
\[ 
\xymatrix{0 \ar[r] & I/I^2 \ar[d]_{\f} \ar[r]^\al & A/I^2 \ar[d]_{g_\f} \ar[r]^\bb & B \ar@{=}[d] \ar[r] & 0\\ 
0 \ar[r] & I/I^2 \ar[r]_{\al_\f} & A_\f \ar[r]_{\bb_\f} & B \ar[r] & 0. } 
\] 
The question is local in $p$. Therefore we can take the completion of \eqref{CD:A_fi} along $p$ to obtain 
\[ 
\xymatrix{0 \ar[r] & \hat B \ar[d]_{\hat \f} \ar[r]^{\al \qquad} & k[[x,y,z]]/(f^2) \ar[d]_{\hat g_\f} \ar[r]^{\qquad\quad\bb} & \hat B \ar@{=}[d] \ar[r] & 0\\ 
0 \ar[r] & \hat B \ar[r]_{\hat \al_{\f}} & \hat{A_\f} \ar[r]_{\hat \bb_{\f}} & \hat B \ar[r] & 0 .}
\] 
Since the construction of $A_\f$ commutes with completion, by 
Corollary \ref{C:embdim} we obtain that $\Hat{A_\f}=A_{\hat \f}$ has 
embedding dimension 3 at $p$ if and only if $\hat \f$ is an 
isomorphism. But completion is faithfully flat, and hence $A_\f$ has 
embedding dimension 3 at $p$ if and only if $\f$ is an isomorphism, 
and so if and only if $g_\f$ is an isomorphism.   
\end{proof}

\subsubsection{Back to the global case} 
 
In this section we want to prove the following:  
\begin{proposition}\label{P:2equiv} 
Let $D$ be a projective 2-dimensional variety over an algebraically 
closed field $k$, having only isolated 
hypersurface-singularities. Suppose $D$ can be embedded into a smooth 
projective 3-fold $X$ and let $\I_D$ denote the ideal 
sheaf of $D$ in $X$. If $H^1(D, \T_D \otimes \II)\simeq H^1(D, 
\mathcal{H}{om}(\Omega_D^1, \II))=0$, then any two embeddings of $D$ 
into smooth three-dimensional varieties, having conormal bundles 
isomorphic to $\II$, are 2-equivalent. 
\end{proposition} 
\begin{proof} 
  
Cover $X$ with affine open sets $\Spec A_i$ such that each of them contains at most one singular point of $D$. Let $V_i := \Spec A_i \cap D=\Spec B_i$ where $B_i=A_i/I_i$.  
 
Consider another embedding $D' \subset X'$ of $D \simeq D'$ into a smooth 3-fold $X'$, with conormal bundle isomorphic to $\I_D/\I_D^2$. 
It gives rise naturally to an infinitesimal extension 
\begin{equation}\label{2D'} 
0 \lright \II \xrightarrow{\al '} \A \xrightarrow{\bb '} \OO_D \lright 0 
\end{equation} 
of $D$ by $\II$, with $\A =\OO_{X'}/\I^2_{D'}$. Then by our previous discussion 
(Proposition \ref{P:embdim3}), for each $i$ there is an isomorphism $\f _i\in \Hom_{B_i}(I_i/I_i^2,I_i/I_i^2)$ such that restricting \eqref{2D'} to $V_i$, we obtain an isomorphism of extensions 
\begin{equation}\label{CD:g_i} 
\xymatrix{0 \ar[r] & I_i/I_i^2 \ar@{=}[d] \ar[r]^{\al_{\f_i}} & A_{\f_i} \ar[d]_{g_i} \ar[r]^{\bb_{\f_i}} & B_i \ar@{=}[d] \ar[r] & 0\\ 
0 \ar[r] & I_i/I_i^2 \ar[r]_{\al'_i} & \A|_{V_i} \ar[r]_{\bb'_i} & B \ar[r] & 0 ,} 
\end{equation} 
where $A_{\f_i}=A_i/I_i^2 \oplus _{\f_i} I_i/I_i^2$ as in \S \ref{ext:general}. 
Let $h_{\f_i}=g_i \circ g_{\f_i}$. Then, by \eqref{CD:A_fi} we obtain the following: 
\[ 
\xymatrix{0 \ar[r] & I_i/I_i^2 \ar[d]_{\f _i} \ar[r]^{\al_i} & A_i/I_i^2 \ar[d]_{h_{\f_i}} \ar[r]^{\bb_i} & B_i \ar@{=}[d] \ar[r] & 0\\ 
0 \ar[r] & I_i/I_i^2 \ar[r]_{{\al'_i}} & \A|_{V_i} \ar[r]_{{\bb'_i}} & B \ar[r] & 0 .} 
\] 
Note that the definition of $h_{\f_i}$ depends on $g_i$. But, according to Lemma \ref{L:trivext}, if we have $g_i$ and $g'_i$ making \eqref{CD:g_i} commutative, then $g'_i=g_i+\al'_i \circ \psi_i \circ \bb_{\f_i}$ with $\psi_i \in  \Der _k(B_i,I_i/I_i^2)$. Hence we get that, correspondingly, $h'_{\f_i}-h_{\f_i}=g'_i \circ g_{\f_i}-g_i \circ g_{\f_i}=(g'_i-g_i)\circ g_{\f_i}=\al'_i \circ \psi_i \circ \bb_i$. 
 
Similarly, suppose that $\f_i$ and $\f'_i$ give isomorphic extensions
$A_{\f _i}$ and $A_{\f'_i}$. Then we know that $\f'_i-\f_i=\psi _i
\circ \al_i$ with $\psi_i \in \Der_k(A_i/I_i^2,I_i/I_i^2)$  (Lemma \ref{L:deriv}), and we have a commutative diagram 
\[ 
\xymatrix{0 \ar[r] & I_i/I_i^2 \ar@{=}[d] \ar[r]^{\al_i} & A_i/I_i^2 \ar[d]_{g_\f} \ar[r]^{\bb_i} & B_i \ar@{=}[d] \ar[r] & 0\\ 
0 \ar[r] & I_i/I_i^2 \ar@{=}[d] \ar[r]^{\al_{\f_i} \circ \f_i} & A_{\f_i} \ar[d]_{\Phi_{\psi_i}} \ar[r]^{\bb_{\f_i}} & B_i \ar@{=}[d] \ar[r] & 0\\ 
0 \ar[r] & I_i/I_i^2 \ar[r]_{\al_{\f'_i}\circ \f'_i} & A_{\f'_i} \ar[r]_{\bb_{\f'_i}} & B_i \ar[r] & 0 } 
\] 
where $\Phi_{\psi_i}$ is defined by $\Phi_{\psi_i}(a,f)=(a,f+\psi_i \circ \al_i \circ \f_i ^{-1} f)$. Then $g_{\f'_i}=\Phi_{\psi_i} \circ g_{\f_i}$ and therefore $h_{\f'_i}=g'_i \circ g_{\f'_i}=(g'_i \circ\Phi_{\psi_i}) \circ g_{\f_i}$ and, as before, differs from  $h_{\f_i}$ by an element of $\Der _k(B_i,I_i/I_i^2)$. 
 
On $V_{ij}=V_i \cap V_j$ we have both $h_{\f_i}$ and $h_{\f_j}$ and by the above 
considerations, we have $h_{\f_i}-h_{\f_j}=\al'_{ij} \circ 
\theta_{ij} \circ \bb_{ij}$ with $\theta_{ij} \in  \Der 
_k(B_{ij},I_{ij}/I_{ij}^2) \simeq \Hom_{B_{ij}} (\Omega _{B_{ij}}, 
I_{ij}/I_{ij}^2)$, and we obtain a well defined cohomology class $[\{\theta_{ij}\}] \in \check{H} ^1 (\{V_i\},  \mathcal{H}{om}(\Omega_D^1, \II)) \simeq H^1(D, \mathcal{H}{om}(\Omega_D^1, \II))$ assigned to the 2-structure $\A$ on $D$. 
 
Conversely, given a cohomology class $[\{\theta  _{ij}\}] \in  \check{H} ^1 (\{V_i\},  \mathcal{H}{om}(\Omega_D^1, \II))$ we can define a ring structure $\A$ on $D$ in the following manner: 
For any $i$, we let  
$\A |_{V_i} := A_i/I_i^2$, where $a_i \in \A|_{V_i}$ is identified with $a_j \in \A|_{V_j}$ if $a_i = a_j +\alpha_{ij} \circ \theta_{ij} \circ \bb_{ij} (a_j)$.  
 
If there is another cocycle $\{\theta'_{ij}\}$ cohomologous with 
$\{\theta_{ij}\}$ and $\A'$ the associated ring structure, then $\A 
\simeq \A'$. Indeed, we have $\theta'_{ij}-\theta_{ij}=\psi_i-\psi_j$ 
with $\psi_i \in \Der_k(B_i,I_i/I_i^2)$, for all $i$, so we can define maps $\Phi_i: A_i/I_i^2 \lright A_i/I_i^2$ by $\Phi_i(a_i)=a_i+\al_i \circ \psi_i \circ \bb_i(a_i)$. Then on $V_{ij}$ we have  
\[ 
\begin{aligned} 
\Phi_i(a_i)&=\Phi_i(a_j +\alpha \circ \theta_{ij} \circ \bb (a_j))=\Phi_i(a_j)+\Phi_i(\alpha \circ \theta_{ij} \circ \bb (a_j))\\ 
&=a_j+\al \circ \psi_i \circ \bb(a_j)+\alpha \circ \theta_{ij} \circ \bb (a_j)+\al \circ \psi_i \circ \bb \circ \alpha \circ \theta_{ij} \circ \bb (a_j)\\ 
&=a_j+\al \circ \psi_j \circ \bb(a_j)+\alpha \circ \theta'_{ij} \circ \bb (a_j)=\Phi_j(a_j)+\alpha \circ \theta'_{ij} \circ \bb (\Phi(a_j)) 
\end{aligned} 
\] 
because $\theta'_{ij}+\psi_j=\theta_{ij}+\psi_i$ and $\bb \circ 
\al=0$. Therefore the $\Phi_i$'s patch to give a morphism of ring 
structures $\Phi :\A \lright \A'$. By a similar construction, we can 
also get  $\Phi' :\A' \lright \A$ so that $\Phi$ and $\Phi'$ will be inverses. 
 
Note that if $[\{\theta  _{ij}\}]=0$, then we get $\A =\OO_X/\I_D^2$,
the 2-structure on $D$.

Therefore, we have the following ``commutative diagram'': 
\[ 
\xymatrix{\{\text{embeddings }D \subset X'\text{ into smooth threefolds}\} \ar[d]\ar[r] & H^1(D, \mathcal{H}{om}(\Omega_D^1, \II)) \ar[dl]\\ 
  \{\text{ring structures locally isomorphic to }\OO_X/\I_D^2\}} 
\] 
where the horizontal arrow is $(D \subset X') \longmapsto \text{the cohomology class 
$[\{\theta_{ij}\}]$}$ defined by \eqref{2D'}, the vertical arrow is $(D 
\subset X') \longmapsto \OO_{X'}/\I^2_{D'}$ and the oblique arrow is 
$[\{\theta_{ij}\}] \longmapsto \A$ as defined above. 
 
Therefore, if $H^1(D, \mathcal{H}{om}(\Omega_D^1, \II)) = H^1(D, \T_D 
\otimes \II)=0$, then any embedding of $D$ into a smooth 3-fold $X'$ 
with ideal sheaf $\I_{D'}$ gives a 2-structure $\OO_{X'}/\I_{D'}^2$ on 
$D$ isomorphic to $\OO_X/\I_D^2$. 
\end{proof}

\section{Reducing the proof of analytic rigidity to the vanishing of 
$H^1(D, {\T}_D)$}\label{red} 
 
We now return to the proof of the Main Theorem.
\begin{theorem}[Main Theorem]\label{T:main1}  Let $X$ be a smooth 3-fold over $\mathbb C$ and suppose we have a $K$-trivial birational 
extremal contraction $\varphi \colon X \longrightarrow~Y$, contracting
a divisor $D$ to a point $q$ in $Y$.  Suppose $D$ is normal, rational,
and $(\omega_D^2) \geq 5$. Then
the contraction $\f$ is analytically rigid.
\end{theorem}

As we showed in \S \ref{known},  under the conditions of the theorem $D \subset X$  is a normal rational del Pezzo surface of 
degree $d \geq 5$ (hence having only $A_n$-type singularities by 
Corollary \ref{C:sing}). We denote by $\pi \colon \DD \lright D$ its 
minimal resolution, and by $\I_D$ the ideal sheaf of $D$ in $X$.
From now on we fix these notations. 
  
Based on the results of \S \ref{first}, the proof of analytic rigidity
is reduced to showing the vanishing of the cohomology groups $H^1(D,
\T_D\otimes \II)$ and $H^1(D, \mathcal T_X \otimes \mathcal I
^{\nu}/\mathcal I ^{\nu +1})$, $\forall \nu \geq 2$. In this section  
we show that $H^1(D, \T_D)=0$ is a sufficient condition for 
achieving these vanishings.
  
\begin{proposition}\label{P:onD}  
If $H^1(D,\T_D)=0$, then 
$H^1(D, \mathcal T_D \otimes \mathcal I_D ^{\nu}/\mathcal I_D ^{\nu 
+1})=0$  for all $\nu \geq 0$.  
\end{proposition}  

\begin{proof}  
First note that $\N_{D/X} \simeq \OO_D(K_D)$ implies  
\begin{equation}\label{E:I/I^2}  
\I_D^\nu/\I_D^{\nu+1} \simeq (\II)^ {\otimes \nu} \simeq \OO_D(- \nu K_D).  
\end{equation}  
  
Let $H \in |-K_D|$ be a general hyperplane section; we can take it to 
be an irreducible, smooth elliptic curve (Proposition \ref{P:cohom}) that 
avoids the singular points of $D$. 
  
Note that $\T_D$ is locally free outside the singular points of $D$, and therefore tensoring
\begin{equation}\label{E:H}  
0 \lright \OO_D(-H) \lright \OO_D \lright \OO_H \lright 0  
\end{equation}  
by $\T_D \otimes \I_D^\nu/\I_D^{\nu+1}$ we obtain
\begin{equation}\label{E:T_D}  
0 \lright \T_D \otimes \I_D^{\nu-1}/\I_D^{\nu} \lright \T_D \otimes \I_D^\nu/\I_D^{\nu+1} \lright \T_D \otimes \I_D^\nu/\I_D^{\nu+1}|_H \lright 0  .
\end{equation}  
  
Next we consider the exact sequence 
\[  
0 \lright \I_H/\I_H^2 \lright \Om|_H \lright \Omega_H^1 \lright 0 , 
\]  
where $\I_H$ is the ideal sheaf of $H$ in $D$.  Since $H$ is an elliptic curve and $\deg(-K_D|_H)=d=\deg(D)$, we
obtain that:
\[  
\begin{aligned}  
h^0&(H,\T_D \otimes \I_D^\nu/\I_D^{\nu+1}|_H)=(2\nu +1)d\\
h^1&(H,\T_D \otimes \I_D^\nu/\I_D^{\nu+1}|_H)=0
\end{aligned}  
\]

Therefore, by induction, using the long exact sequence associated to 
\eqref{E:T_D} and the condition $H^1(D,\T_D)=0$, we obtain $H^1(D, 
\mathcal T_D \otimes \mathcal I_D ^{\nu}/\mathcal I_D ^{\nu +1})=0$ 
for all $\nu \geq 0$.  
\end{proof}

\begin{corollary}\label{CC:onX} 
Under the assumptions of the Main Theorem, if $H^1(D,\T_D)=0$, then any two embeddings of $D$ 
into smooth three-dimensional varieties, having conormal bundles 
isomorphic to $\II$, are 2-equivalent
\end{corollary} 
 
\begin{theorem}\label{P:onX} 
Under the assumptions of the Main Theorem, if $H^1(D,\T_D)=0$, then $H^1(D, \mathcal T_X \otimes \mathcal I_D ^{\nu}/\mathcal I_D ^{\nu +1})=0$, $\forall \nu \geq 2$.  
\end{theorem}  
\begin{proof} 
In the case when $D$ is nonsingular, the teorem is an easy consequence of
Proposition \ref{P:onD} and the vanishing of the first cohomology of
$\OO _D ((1-\nu) K_D)$ for any del Pezzo surface  (Proposition \ref{P:cohom}). We
actually obtain the vanishing of $H^1(D, \mathcal T_X \otimes \mathcal
I_D ^{\nu}/\mathcal I_D ^{\nu +1})=0$, for all $\nu \geq 1$.  

So we concentrate on the case when $D$ is normal rational. In this
case the difficulty comes from the fact that the tangent sheaf 
$\T_D$ is not locally free anymore. In fact, 
$H^1(D,\T_X \otimes \II)$ doesn't vanish in general (see Remark 
\ref{R:I/I^2} below), but an easy argument shows that 
it is enough to prove that $H^1(D,\T_X \otimes \I_D^2/\I_D^3)=0$.  

Indeed, for a general $H \in |-K_D|$, as in the proof of Proposition
\ref{P:onD}, tensor
\[  
0 \lright \OO_D(-H) \lright \OO_D \lright \OO_H \lright 0  
\]  
by $\T_X \otimes \I_D^\nu/\I_D^{\nu+1}$ for $\nu \geq 2$. From the long exact sequence of cohomology, we get 
\[  
H^1(D, \T_X \otimes \I_D^{\nu -1}/\I_D^{\nu}) \lright H^1(D, \T_X \otimes \I_D^\nu/\I_D^{\nu+1}) \lright H^1(H, \T_X \otimes \I_D^\nu/\I_D^{\nu+1}|_H).  
\]  
  
Consider the exact sequences  
\[  
0 \lright \T_H \lright \T_X|_H \lright \N_{H/X} \lright 0  
\]  
and  
\[  
0 \lright \N_{H/D} \lright \N_{H/X} \lright \N_{D/X}|_H \lright 0 . 
\]  
The latter sequence is exact because both $H$ in $D$ and $D$ in $X$ 
are Cartier divisors; see for example \cite[19.1.5]{Groth2}. 
 
We know that $\N_{H/D}=\OO _H(-K_D)$, $\N_{D/X}|_H=\OO _H(K_D)$, so 
tensoring the above exact sequences with 
$\I_D^\nu/\I_D^{\nu+1}=\OO_D(-\nu K_D)$ yields $H^1(H,\N_{H/X} \otimes 
\I_D^\nu/\I_D^{\nu+1})=0$ (because $-K_D$ is nonzero and effective) 
and hence $H^1(H, \T_X \otimes \I_D^\nu/\I_D^{\nu+1}|_H)=0$ for all $\nu 
\geq 2$.  (Recall that $H$ is an elliptic curve.)
  
Therefore we have a surjection  
\[  
H^1(D, \T_X \otimes \I_D^{\nu -1}/\I_D^{\nu}) \lright H^1(D, \T_X 
\otimes \I_D^\nu/\I_D^{\nu+1}) \lright 0 \qquad \forall \nu \geq 2,
\]  
and so we can use induction on $\nu$ to prove the vanishing of $H^1(D, \T_X \otimes \I_D^\nu/\I_D^{\nu+1})$.   

\noindent The proof is finished by the following proposition.  
\end{proof}  
  
\begin{proposition}\label{L:onX}  
Suppose $H^1(D,\T_D)=0$. Then $H^1(D,\T_X \otimes \I_D^2/\I_D^3)=0$.  
\end{proposition}  
\begin{proof}  
The sequence  
\begin{equation}\label{Eq:conormal}  
0 \lright \II \lright \Omega^1_X|_D \lright \Om  \lright 0  
\end{equation}  
is exact because $D$ is a Cartier divisor in the smooth 3-fold 
$X$. Applying the functor $\mathcal Hom_{\OO_D}( -, \I_D^2/\I_D^3)$ we obtain 
\[  
\begin{aligned}  
0 \lright &\mathcal Hom(\Om, \I_D^2/\I_D^3) \lright \mathcal Hom(\Omega^1_X|_D, \I_D^2/\I_D^3) \lright \mathcal Hom(\II, \I_D^2/\I_D^3) \lright \\  
\lright &\mathcal Ext^1(\Om , \I_D^2/\I_D^3) \lright \Ex(\Omega^1_X|_D, \I_D^2/\I_D^3)=0,  
\end{aligned}  
\]  
and hence  
\[  
\begin{aligned}  
0 \lright &\T_D \otimes \I_D^2/\I_D^3 \lright \T_X|_D \otimes \I_D^2/\I_D^3 \lright \mathcal Hom(\II, \I_D^2/\I_D^3)  
\lright &\mathcal Ext^1(\Om , \I_D^2/\I_D^3) \lright 0 . 
\end{aligned}  
\]  
Break this up into two short exact sequences  
\begin{equation}\label{Eq:1}  
0 \lright \T_D \otimes \I_D^2/\I_D^3 \lright \T_X \otimes \I_D^2/\I_D^3 \lright \mathcal C \lright 0  
\end{equation}  
and  
\begin{equation}\label{Eq:2}  
0 \lright \mathcal C \lright  \mathcal Hom(\II, \I_D^2/\I_D^3) \lright \mathcal Ext^1(\Om , \I_D^2/\I_D^3) \lright 0 . 
\end{equation}  
By Proposition \ref{P:onD}, $H^1(D,\T_D)=0$ implies $H^1(D,\T_D 
\otimes \I_D^2/\I_D^3)=0$. (This is the only place in the proof where we 
are using the assumption $H^1(D,\T_D)=0$.) As $\mathcal Hom(\II,
\I_D^2/\I_D^3) \simeq \II \simeq \OO_D(-K_D)$, we also have $H^1(D,\mathcal Hom (\II, 
\I_D^2/\I_D^3)) =  H^1(D,\OO_D(-K_D))=0$ (\cite{Wat}). Hence, from the corresponding long exact sequences on cohomology we get 
\[  
0 \lright H^1(D,\T_X \otimes \I_D^2/\I_D^3) \lright H^1(D,\mathcal C) = 
\coker(\Phi),  
\]  
where $\Phi :\Hom (\II, \I_D^2/\I_D^3) \lright H^0(D, \mathcal Ext^1(\Om ,  
\I_D^2/\I_D^3))$ is obtained from 
\eqref{Eq:2} by taking global sections.
  
So now our goal is to show that $\Phi$ is surjective, and hence obtain
that $H^1(D,\T_X \otimes \I_D^2/\I_D^3)=0$. As $\dim H^0(D,\Ex(\Om ,
\I_D^2/\I_D^3)) =\Sigma_{i=1}^{r} \lambda_i$ if $D$ has $r$
singularities, of type $A_{\lambda_1}$, $A_{\lambda_2}$, ...,
$A_{\lambda_r}$ respectively (Lemma \ref{L:Ext}), and $\Hom (\II, \I_D^2/\I_D^3) \simeq
H^0(D,\OO_D(-K_D))$, it is enough to show that $\dim \ker(\Phi) =\dim (H^0(D,\OO_D(-K_D)))
-\Sigma_{i=1}^{r} \lambda_i$.

Let us 
look  at the kernel of $\Phi$. Retaining the notations from 
\S \ref{2str}, we observe (using \v{C}ech cohomology), that $\Phi$ is given by  
\[  
s \mapsto \Big\{(0 \lright \I_D^2/\I_D^3 \lright \mathcal M_s \lright \Om  \lright 0)|_U\Big\}_U \in \check H^0(\U,\mathcal Ext^1(\Om , \I_D^2/\I_D^3))  
\]  
for any morphism $s \in \Hom (\II, \I_D^2/\I_D^3)$, where $\mathcal M_s =\Omega^1_X|_D \oplus _s \I_D^2/\I_D^3$ from the following diagram:  
\begin{equation}\label{CD:M_s}  
\xymatrix{0 \ar[r] & \II \ar[d]_{s} \ar[r] & \Omega^1_X|_D \ar[d] \ar[r] & \Om \ar@{=}[d] \ar[r] & 0\\  
0 \ar[r] & \I_D^2/\I_D^3  \ar[r] & \mathcal M_s \ar[r] & \Om \ar[r] & 0 .}  
\end{equation}   
Therefore $\Phi(s)=0$ is a local condition, equivalent to the exact sequence $0 \lright \I_D^2/\I_D^3 \lright \mathcal M_s 
\lright \Om  \lright 0$ splitting locally. On the smooth portion of $D$, 
$\mathcal Ext^1(\Om , \I_D^2/\I_D^3)$ vanishes, since $\Om$ is locally 
free.  
Hence we only have to look near the 
singular points of $D$. Let $p$ be a singular point of $D$, and let $V=\Spec B \subset D$ be an affine open neighborhood of $p$. 
 
Locally then $s \in \Hom (\II, \I_D^2/\I_D^3)$ is a homomorphism $s:B \lright B$, and the diagram \eqref{CD:M_s} translates to  
\[  
\xymatrix{0 \ar[r] & B \ar[d]_{s} \ar[r]^{\alpha} & B^3 \ar[d]_{\gamma} \ar[r]^{\beta} & M \ar@{=}[d] \ar[r] & 0\\  
0 \ar[r] & B \ar[r]_{\alpha_s \quad} & B^3\oplus _s B \ar[r]_{\quad \beta_s} & M \ar[r] & 0 }  
\]  
where $M=\Om|_V \simeq \Omega_B$. The bottom row splits if and only if there 
is a map $P: B^3\oplus _s B \lright B$ such that $P \circ 
\alpha_s=id_B$, i.e. $P(0,g)=g$, for all $g \in B$. This is equivalent 
to the existence of $Q: B^3 \lright B$ such that $Q \circ \alpha =s$.  
  
Indeed, if $P$ exists, then take $Q:=P \circ \gamma$. Conversely, given $Q$, define $P(v,g):=Q(v)+g$.  
  
Therefore the bottom row splits if and only if $s \in \im\alpha ^*$ in  
\[  
0 \lright \Hom(M,B) \xrightarrow{\beta ^*} \Hom(B^3,B) \xrightarrow{\alpha  ^*} \Hom(B,B).   
\]  
  
Completion is faithfully flat, so we can replace $B$ by its completion 
at the maximal ideal of $B$ corresponding to $p$. Recall that $D$ has
only $A_\lambda$-type singularities (Corollary \ref{C:sing}), hence we can assume 
that $B =k[[x,y,z]]/(f)$ where $f(x,y,z)=x^{\lambda+1}+yz$ 
corresponding to an $A_\lambda$-singularity . In this case, the map 
$\al \colon B \lright B^3$ is given by the partial derivatives of 
$f$: $\al(1)=\left(\begin{smallmatrix}  
(\lambda+1)x^\lambda\\  
y\\  
z  
\end{smallmatrix}\right)$. Therefore, if we identify $\Hom(B^l,B)$ 
with $B^l$, the condition $s \in \im\alpha ^*$  is equivalent to $s \in 
(x^\lambda,y,z)B$, and hence to $s \in J(p)$, where  $J(p)$ is the Jacobian ideal of $D$ at $p$. Note that this is independent of a coordinate change.    
  
Therefore, $\ker(\Phi)=\{s \in \Hom (\II, \I_D^2/\I_D^3)\simeq
H^0(D,\OO_D(-K_D))$: $s \in J(p)$ for
any singular point $p$ of $D$\}, and the proof is concluded by the following:

\begin{claim} If $D$ has $r$ singularities, of type $A_{\lambda_1}$,
$A_{\lambda_2}$, ..., $A_{\lambda_r}$, then, using the above
description of   $\ker(\Phi)$,  
$\dim \ker(\Phi) =\dim H^0(D,\OO_D(-K_D)) -\Sigma_{i=1}^{r} \lambda_i$.  
\end{claim}  
It is easy to see the Claim is true in the case when $D$ has one or two $A_1$-singularities, 
because in this case $s$ belongs to the Jacobian ideal if and only if the corresponding hyperplane section passes through the (two) 
singular point(s). However, the statement is not obvious for any other 
combination of singularities for $D$, and the proof will be given in Lemma 
\ref{Claim:onX} below. 
\end{proof}  
  
\begin{lemma}\label{L:Ext}  
Let $D$ be a projective surface over an algebraically closed field $k$ 
of characteristic 0. Suppose $D$ has $r$ singularities, of types $A_{\lambda_1}$, 
$A_{\lambda_2}$, ..., $A_{\lambda_r}$ respectively. Then $h^0(D,\Ex(\Om , \I_D^2/\I_D^3)) =\Sigma_{i=1}^{r} \lambda_i$.  
\end{lemma}  
\begin{proof}  
Since $\Ex  (\Om , \I_D^2/\I_D^3))$ is a skyscraper sheaf, we may assume that $D$ is affine, say $D=\Spec B$, having one singularity $p$, of type $A_{\lambda}$. Then the completion $\hat B_p$ of the local ring $B_p$ is isomorphic to $k[[x,y,z]]/(x^{\lambda+1}+y^2+z^2)$. Consider the exact sequence:  
\[  
0 \lright \hat B_p \xrightarrow{\left(\begin{smallmatrix}  
(\lambda+1)x^\lambda\\  
y\\  
z  
\end{smallmatrix}\right)} \hat B_p ^3 \lright \Omega_{\hat B_{p/k}} \lright 0.  
\]  
From here, using the long exact sequence on Ext's, we can easily conclude that  $\Ext^1_{\hat B_{p}} 
(\Omega_{\hat B_p/k}, \hat B_p) \simeq 
(k[[x,y,z]]/(x^{\lambda+1}+y^2+z^2))/(x^{\lambda},y,z) 
\simeq k^{\lambda}$. But $\Omega_{\hat B_p/k} \simeq \Omega_{B_p/k} \otimes 
_{B_p} \hat B_p$ and completion is faithfully flat, and  therefore 
we obtain an isomorphism $\Ex  (\Om ,\OO_D) \simeq k^{\lambda}$.   
\end{proof}  

\begin{lemma}\label{Claim:onX} 
Let $D$ be a normal, rational del Pezzo surface of degree $d \geq 5$. If $D$ has $r$ singularities, $p_1, p_2, \dots p_r$, of type
$A_{\lambda_1}$, $A_{\lambda_2}$, ..., $A_{\lambda_r}$ respectively,
let us denote by $J(p_i)$ the Jacobian ideal of $D$ at $p_i$. Then  
$$\dim \{s \in H^0(D,\OO_D(-K_D)):\, s \in J(p_i) \, ,\forall i\}=\dim
H^0(D,\OO_D(-K_D)) -\Sigma_{i=1}^{r} \lambda_i .$$
\end{lemma}
 \begin{proof}
As observed above, the lemma is obvious in the case when $D$ has one or two $A_1$-singularities, 
because in this case $s$ belongs to the Jacobian ideal if and only if
the corresponding hyperplane section passes through the (two) singular
point(s). So we concentrate on the cases when $D$ is more singular. 

In order to compute $\dim  \{s \in H^0(D,\OO_D(-K_D)):\, s \in
J(p_i) \, ,\forall i\}$, we look at the anticanonical embedding of $D$
into $\PP^d$ (\cite[Cor. 4.5]{Wat}). An element $s$ of
$H^0(D,\OO_D(-K_D))$ can be considered as a hyperplane section of $D
\subset \PP^d$.

If the minimal resolution $\DD$ of $D$ is obtained from  
the blowing up of a set of (possibly infinitely near) points $\Sigma$ on $\PP^2$
(as described in \S \ref{delp}), we have $H^0(D, \OO_D(-K_D)) \simeq H^0(\DD, \OO_\DD(-K_\DD))=\big\{s \in  
H^0(\PP^2,\OO_{\PP^2}(3)) \colon s \text{ has}\\\text {base-points at the points  
of }\Sigma\big\}$. Therefore the anticanonical embedding  of $D$ is the 
closure of the image of a map $\f_0 \colon \PP^2 \cdots \rightarrow \PP^d$ having  
the assigned base-points $\Sigma$. Fixing homogeneous coordinates $[x_0:x_1:x_2]$  
on $\PP^2$ and $[a:b:c:d:e:\dots]$ on $\PP^d$, this map can be explicitly  
described, and we can find the homogeneous ideal $I_0$ of the closure of its image by  
elimination theory (using the computer programs Maple or Macaulay2).
Therefore we have an explicit description of the homogeneous ideal defining $D
\subset \PP^d$, and so we can compute the complete local rings of $D$ at the
singular points of $D$ and hence their respective Jacobian ideals. 

In the following we describe in detail the case when the minimal resolution $\DD$ 
of $D$ is isomorphic to $S_4''$ (notation as in \S \ref{delp}). For the rest of the  
cases, the computations are similar, and are left to the reader.

The degree of $D$ when  $\DD=S_4''$  is $d=5$, therefore $|-K_D|$ embeds 
$D$ into $\PP ^5$ via $\f_0 \colon \PP^2 \cdots \rightarrow \PP^5$.
By the description of the blow-up from \S \ref{delp}, $ H^0(D,\OO_D(-K_{D}))$ is obtained from those sections of  
$H^0(\PP^2,\OO_{\PP^2}(-K_{\PP ^2}))$ that pass through the infinitely  
near points given by $\Sigma_4''$, namely $[x_0:x_1:x_2]=[1:0:0],
[1:0:1]$, $[1:0:-1]$, and the point $p_4$ corresponding to the tangent direction
to the line $x_2=0$, i.e. from those sections $s \in H^0(\PP^2,\OO_{\PP^2}(-K_{\PP ^2}))$ for which we have:  
\[\begin{aligned}  
s(1,0,0) &=0\\  
s(1,0,1) &=0\\  
s(1,0,-1) &=0\\  
s_{x_1}(1,0,0)&=0.  
\end{aligned}  \]  
These give us the map $\varphi _0 : \PP ^2 \cdots \rightarrow \PP ^5$,
defined by   
\[  
[x_0:x_1:x_2] \lright [x_1^3: x_2^3-x_0^2x_2: x_0x_1^2: x_1^2x_2:
x_1x_2^2: x_0x_1x_2].  
\]   
The image of  $\varphi_{|-K_D|}$ in $\PP ^5$ is the closure of  $\varphi_0(\PP ^2)$.

To compute the ideal of $D \simeq \overline{\varphi_0(\PP ^2)} \subset
\PP^5$, we use
elimination, with the help of the Maple or Macaulay 2 computer
programss. We obtain that the embedding of $D$ into $\PP^5$ is given by the
ideal
\[ 
I_0=(ce - df, d^2 - ae, cd - af, bd - e^2 + f^2, ab - de + cf),
\] 
where the homogeneous coordinates in $\PP ^5$ are  given by
$a,b,c,d,e,f$. In this case $D$ has one $A_2$-singularity, at $p=[0:1:0:0:0:0]$.  

The singular point of $D$ is in the affine piece $\{b \neq 0 \}$ of
$\PP ^5$. After dehomogenizing (putting $b=1$ and keeping the
notation unchanged for the other coordinates), we obtain that the coordinate ring of
$D$ at $p$ is $B=k[c,e,f]/(ce-fe^2+f^3)$ (where $a=-cf+e^3-f^2e$ and $d=e^2-f^2$). Therefore the Jacobian ideal
at $p$ is $J(p)=(c,e,f^2)$. Direct computation shows that $\{s \in H^0(D,\OO_D(-K_D)):\, s \in J(p) \, \}=\{\mu_1 a+ \mu_2 b+ \mu_3 
c+ \mu_4 d+ \mu_5 e+ \mu_6 f : \mu_2=\mu_6=0\}$, hence it has the
required dimension.
\end{proof} 
\begin{corollary}\label{C:main} 
Suppose that $H^1(D,\T_D)=0$. Then the Main Theorem holds, i.e. the contraction $\f$ is analytically rigid.  
\end{corollary}  
 
\begin{proof}
Suppose $D$ and $X$ are as in the Main Theorem (Theorem \ref{T:main1}). Then, if $H^1(D, 
\T_D)=0$,
any two embeddings of $D$ into smooth 3-folds with normal bundles 
isomorphic to $\OO_D(K_D)$ are 2-equivalent and we also have the 
vanishing of $H^1(D, \mathcal T_X \otimes \mathcal I ^{\nu}/\mathcal I 
^{\nu +1})$ for any $\nu \geq 2$ (Corollary \ref{CC:onX} and Theorem \ref{P:onX}). Therefore, 
by Lemma \ref{L:Hiron} and \cite[Theorem 3]{Hir}, the Main Theorem holds
\end{proof}
\begin{remark}\label{R:I/I^2} In the case when $D$ is nonsingular, the
vanishing $H^1(D,\T_X 
\otimes \II)=0$ is and easy consequence of the vanishing of $H^1(D, 
\mathcal T_D)$ (via the exact sequence \eqref{Eq:conormal}, using 
Proposition \ref{P:onD} and the vanishing of $H^1(D,\OO_D)$). Hence,
by induction (as in the proof of Theorem \ref{P:onX}), 
we have the vanishing of $H^1(D, \mathcal T_X \otimes \mathcal I_D 
^{\nu}/\mathcal I_D ^{\nu +1})$, for all $\nu \geq 2$. However, if $D$ is singular, it is not generally true that $H^1(D, 
\mathcal T_D)=0$ implies $H^1(D,\T_X \otimes \II)=0$.  
\end{remark}  
  
Indeed, suppose $D$ has at least one singular point. Consider the exact sequence   
\[  
\begin{aligned}  
0 \lright &\Hom(\Om,\II) \lright \Hom(\Omega^1_X|_D, \II) \lright \Hom (\II,\II) \xrightarrow{\delta} \\  
\xrightarrow{\delta} &\Ext ^1(\Om,\II) \lright \Ext ^1(\Omega^1_X|_D, \II) \lright \Ext ^1(\II,\II)  
\end{aligned}  
\]  
obtained from \eqref{Eq:conormal}. By Serre duality and $D$ being del 
Pezzo, we have the following: $\Ext ^1(\II,\II) \simeq
H^1(D,\OO_D)=0$, $\Ext ^1(\Omega^1_X|_D, \II) \simeq H^1(D,\T_X|_D \otimes \II)$ and $\Hom (\II,\II) \simeq k$, so we obtain  
\[  
\begin{aligned}  
0 \lright &\Hom(\Om,\II) \lright \Hom(\Omega^1_X|_D, \II) \lright k \xrightarrow{\delta} \\  
\xrightarrow{\delta} &\Ext ^1(\Om,\II) \lright H^1(D,\T_X|_D \otimes \II) \lright 0.  
\end{aligned}  
\]  
  
But $\delta (1)$ corresponds to the extension \eqref{Eq:conormal} that is not split, and hence not zero in $\Ext ^1(\Om,\II)$. Therefore $\delta$ is injective and we have  
\begin{equation}\label{Eq:A_1}  
0 \lright  k \xrightarrow{\delta} \Ext ^1(\Om,\II) \lright H^1(D,\T_X \otimes \II) \lright 0.  
\end{equation}  
  
From the five-term exact sequence associated to the local to global spectral sequence, we have  
\[\begin{aligned}  
0 \lright &H^1(D, \mathcal{H}{om}(\Omega_D^1, \OO_D)) \lright \Ext ^1  
 (\Omega_D^1,\OO_D ) \lright \Gamma(D,\Ex  (\Omega_D^1,\OO_D ))\lright \\ \lright  
&H^2(D, \mathcal{H}{om}(\Omega_D^1,\OO_D )) \lright \Ext ^2 (\Omega_D^1, \OO_D)  
\end{aligned}\]  
or with the notation $\T_D=\mathcal Hom (\Om, \OO_D)$ and using the vanishing 
$H^2(D,\T_D)=0$ (\cite[Lemma 5.6]{Gr}):  
\begin{equation}\label{E:Local-global}  
0 \lright H^1(D,\T_D) \lright \Ext ^1 (\Omega_D^1,\OO_D ) \lright \Gamma(D,\Ex  (\Omega_D^1,\OO_D )) \lright 0.  
\end{equation}  
  
We know by Lemma \ref{L:Ext} that $\dim \Gamma(D,\Ex 
(\Omega_D^1,\OO_D ))=\Sigma_{i=1}^{r} \lambda_i$, if $D$ has $r$ 
singularities, of type $A_{\lambda_1}$, $A_{\lambda_2}$, ..., 
$A_{\lambda_r}$ respectively. Therefore, if $H^1(D, \T_D)=0$, we obtain 
$\dim \Ext ^1 (\Omega_D^1,\OO_D )=\Sigma_{i=1}^{r} \lambda_i$.  
 
Now 
suppose that 
$H^1(D,\T_X \otimes \II)=0$. Then the exact sequence \eqref{Eq:A_1} implies that $\dim 
\Ext ^1 (\Omega_D^1,\OO_D )=1$, and hence $D$ can have only one 
singularity, of type $A_1$. However, this is not the case in general, 
and hence $H^1(D,\T_X \otimes \II)$ doesn't vanish for a general singular del Pezzo surface with  $H^1(D,\mathcal T_D)=0$.

\section{Computing the obstruction to formal equivalence: the 
vanishing of $H^1(D,{\mathcal T}_D)$ 
}\label{vanish}

Let $D$ be a normal rational del Pezzo surface of degree $d \geq 5$.  We showed (Corollary \ref{C:main}) that 
$H^1(D,\T_D)=0$ is a sufficient condition for the analytic rigidity of 
the contraction $\f$ in the Main Theorem. Here we show the vanishing of $H^1(D,\T_D)$.

\subsection{The Leray spectral sequence}\label{Leray}

As $D$ is
a surface having only isolated normal singularities, 
$\pi_* \T_\DD \simeq \T_D$ (\cite[Prop. (1.2)]{B-W} ). We use this fact and the Leray 
spectral sequence $$E_2 ^{p,q} =H^p (D, R^q \pi_* \T_\DD) 
\Rightarrow E_{\infty}^{p+q}=H^{p+q}(\tilde D,\T_\DD)$$ 
 to compare the cohomology of $\T_D$ with that of $\T_\DD$.   
  
The first four terms of  the corresponding five-term exact sequence (\cite[Theorem I.4.5.1]{Godement})  
\[  
0 \lright E_2 ^{1,0} \lright E_{\infty}^1 \lright E_2 ^{0,1} \lright 
E_2 ^{2,0}\lright  E_{\infty}^2
\]  
give in our case  
\[  
0 \lright H^1(D, \pi_* \T_\DD) \lright H^1 (\tilde D,\T_\DD) \lright H^0 (D, R^1 \pi_* \T_\DD) \lright H^2(D, \pi_* \T_\DD). 
\]  
With the identification $\pi_* \T_\DD \simeq \T_D$, and using 
$H^2(D,\T_D)=0$ \cite[Lemma 5.6]{Gr}, we obtain
\begin{equation}\label{Ex:DD} 
0 \lright H^1(D, \T_D) \lright H^1 (\tilde D,\T_\DD) \lright H^0 (D, R^1 \pi_* \T_\DD) \lright 0 . 
\end{equation} 

In order to show the vanishing of $H^1(D, \T_D) $, we  show that
$\dim H^1 (\tilde D,\T_\DD) = \dim H^0 (D, R^1 \pi_* \T_\DD)$.  
 
\subsection{Local computations: $H^0(D,R^1 \pi_* \T_\DD)$}\label{local} 
 As $D$ is a normal rational del Pezzo surface of degree $d \geq 
 5$, it has only singularities of type $A_{1}$, $A_{2}$, $A_{3}$ 
 and $A_{4}$ (Corollary \ref{C:sing}).  Denote by $E$ the exceptional
locus of  its minimal resolution $\pi \colon \DD \lright D$.  
  
We have $H^0(D,R^1 \pi_* \T_\DD)\simeq R^1 \pi_* \T_\DD$  if regarded
as a complex vector-spaces, as $R^1 \pi_* (\pi ^* \T_\DD)$ is a
skyscraper sheaf supported on the singular points of $D$.  
  By the theorem of formal functions (\cite[Theorem III.11.1]{Hartsh}), we obtain   
\[  
R^1 \pi_* \mathcal \T_\DD \simeq \varprojlim H^1(\mathcal E_n, \mathcal \T_\DD|_{\EE_n}),  
\]  
where $\mathcal E_n$ is the closed subscheme of $\tilde D$ defined by $\mathcal I^n_E$, where $\mathcal I_E$ is the ideal sheaf of $E$ in $\tilde D$.  
  
From \cite[(1.6)]{B-W}, we have the following lemma, true for any surface
having only isolated rational singularities: 
\begin{lemma}\label{L:B-W} 
If $Z$ is an effective divisor on $\DD$ supported on $E$, there is an exact sequence  
\begin{equation} \label{Ex:B-W} 
0 \lright \T_Z \lright \T_\DD|_Z \lright \bigoplus_{i=1}^\lambda \N_{E_i/\DD} \lright 0  
\end{equation}  
where $E_1, E_2, \dots ,E_\lambda$ are the irreducible components of 
$E$ and $\N_{E_i/\DD} :=\OO_{E_i}(E_i)=\OO_{E_i}(-2)$ is the normal 
bundle of $E_i$ in $\DD$.  (The second map of \eqref{Ex:B-W} is the sum of the compositions $\T_\DD \otimes \OO_Z 
\lright \T_\DD \otimes \OO_{E_i} \lright \N_{E_i/\DD}.$)
\end{lemma}  

By the tautness of rational double point singularities (\cite{Tju}) we
have that $H^1(\T_Z)=0$. 
Therefore, the long exact sequence obtained from \eqref{Ex:B-W}
implies $h^1(\mathcal E_n, \mathcal \T_\DD|_{\EE_n})=h^1(\mathcal E_n,
\bigoplus_{i=1}^\lambda \N_{E_i})$, for all $n \geq 0$. This shows that
if $D$ has $r$ singularities, of type $A_{\lambda_1}$, $A_{\lambda_2}$, 
..., $A_{\lambda_r}$ respectively, then $h^0(D,R^1 \pi_* \T_\DD)=\Sigma_{i=1}^{r} \lambda_i$.

\subsection{Global computations: $H^1 (\tilde D,\T_\DD)$} Here we show that $\dim 
H^1(\DD,\T_\DD)=\Sigma_{i=1}^{r} \lambda_i$ as well, and therefore we
obtain $H^1(D,\T_D)=0$ from \eqref{Ex:DD} .

First, we need some preliminary results relating the 
tangent bundle of a (smooth) surface $S$ to that of a one-point 
blow-up of $S$ (Lemma \ref{P:coker}). 
  
\begin{lemma}\label{P:blup}  
Let $\s \colon S' \lright S$ be a birational morphism of smooth projective surfaces and let $\F$ be a locally free sheaf on $S$. Then $H^*(S, \F)=H^*(S', \s^* \F)$.  
\end{lemma}  
\begin{proof}  
The morphism $\s$ can be factored as the composition of blow-ups. Therefore it is enough to assume $\s$ itself is the blow-up of a point $p \in S$.  
  
Because the sheaf $\F$ is locally free, the projection formula and 
normality of $S$ imply that $\s_* \s^* \F \simeq \F$. Therefore, by 
a degenerate case of the Leray spectral sequence, it is enough to show 
that $R^i \s_* \s^* \F=0$, $\forall i >0$. But again, by the 
projection formula and $\s_* \OO_{S'} =\OO_S$, we can reduce this to
showing that $R^i \s_* \OO_{S'}=0$, which is proven in \cite[Prop. V.3.4]{Hartsh}. 
\end{proof}  
  
\begin{lemma}\label{P:coker}  
Let $\s \colon S' \lright S$ be the blowing up of a smooth projective 
surface $S$ at a point $p$  and let $\EE$ denote the exceptional locus. We then have an exact sequence  
\begin{equation}\label{blup1}  
0 \longrightarrow \T _{S'}  \longrightarrow \sigma ^* \T_{S} \longrightarrow \OO_{\EE}(1) \longrightarrow 0.  
\end{equation}  
\end{lemma}  
\begin{proof}  
By smoothness of $S$ at $p$, we may assume (after completion) that $p$
 has a neighborhood in $S$ analytically  isomorphic to $\mathbb
A^2$. In this case we may consider $\s$ as  the blow-up $\s \colon U \lright \mathbb A^2$ of the
origin in $\mathbb A^2$. The lemma now follows from an easy
computation in local coordinates.
\end{proof}

Here is also a more general version:  
\begin{lemma}\label{P:coker1}  
Let $\s \colon S' \lright S$ be the blowing up of a smooth projective 
variety $S$ of dimension $n$ at a point $p$  and let $\EE$ denote the exceptional locus. We then have an exact sequence  
\begin{equation}
0 \longrightarrow \T _{S'}  \longrightarrow \sigma ^* \T_{S}
\longrightarrow \T_\EE \otimes \OO_{\EE}(\EE) \longrightarrow 0.  
\end{equation}  
\end{lemma}  
\begin{proof}  
As before, we may assume (after completion) that $p$
 has a neighborhood in $S$ analytically  isomorphic to $\mathbb
A^n$, and that $\s \colon U \lright \mathbb A^n$ is the blow-up of the
origin in $\mathbb A^n$.
  
We have the first fundamental exact sequence of differentials:
\begin{equation}\label{E:diff} 
0 \lright \s^* \Omega_{\mathbb A^n}^1 \lright \Omega_U^1 \lright \Omega_{U/{\mathbb A^n}} \lright 0  
\end{equation} 
which gives, after taking $\mathcal Hom_{\OO_U}(-,\OO_U)$: 
\[ 
0 \longrightarrow \T _U  \longrightarrow \sigma ^* \T_{\mathbb A ^n} \longrightarrow \EE xt^1_{\OO_U}(\Omega^1_\EE,\OO_U) \lright 0,  
\]  
where we used $\Omega_{U/{\mathbb A^n}} \simeq \Omega^1_\EE$ (see for example \cite{Kleinman}). 
  
We can compute $\EE xt^1_{\OO_U}(\Omega^1_\EE,\OO_U)$ using the conormal exact sequence  
\begin{equation}\label{E:conormal}  
0 \lright \I_\EE/\I_\EE^2 \lright \Omega_U^1|\EE \lright \Omega^1_\EE \lright 0.  
\end{equation}  
The sheaf $\I_\EE/\I_\EE^2$ is supported on $\EE$, hence $\mathcal 
Hom_{\OO_U}(\I_\EE/\I_\EE^2,\OO_U)=0$. Also, we have $\EE 
xt^2_{\OO_U}(\Omega^1_\EE,\OO_U)=0$ because \eqref{E:diff} is a 
locally free resolution of $\Omega^1_\EE$ of length 1. Therefore, 
using $\I_\EE \simeq \OO_U(-\EE)$, we obtain:  
\begin{equation}\label{blup2}  
0 \lright \EE xt^1_{\OO_U}(\Omega^1_\EE,\OO_U) \lright \EE xt^1_{\OO_U}(\Omega^1_U|_\EE,\OO_U) \lright \EE xt^1_{\OO_U}(\OO_U(-\EE)|_\EE,\OO_U) \lright 0.  
\end{equation}

From the exact sequence 
\[  
0 \lright \OO_U(-\EE) \lright \OO_U \lright \OO_\EE \lright 0  
\]  
we conclude that $\EE xt^1_{\OO_U}(\OO_\EE,\OO_U)=\OO_\EE(\EE)$, and hence, from \eqref{blup2}, we have   
\[  
0 \lright \EE xt^1_{\OO_U}(\Omega^1_\EE,\OO_U) \lright \T _U \otimes _{\OO_U}\OO_\EE(\EE) \lright \OO_\EE(2\EE) \lright 0.  
\]  
Tensoring the dual of \eqref{E:conormal} by $\OO_\EE(\EE)$, we obtain
\[  
0 \lright \T_\EE \otimes  \OO_\EE(\EE) \lright \T _U \otimes _{\OO_U}\OO_\EE(\EE) \lright \OO_\EE(2\EE) \lright 0.  
\]  
Comparing the last two exact sequences, we obtain the desired result.
\end{proof}

Lemma \ref{P:coker} gives us a tool to compute $H^1 (\tilde 
D,\T_\DD)$ step-by-step, blowing up one point at a time: 
 
Let $S=V(\Sigma)$ be a surface obtained by blowing up a set $\Sigma$ of (possibly infinitely near) points 
in almost general position 
on $\PP^2$. 
Let $p$ be a point on $S$ such that the 
points of the  set $\Sigma \cup \{p\}$ are again in almost general 
position, and denote $S'=V(\Sigma \cup \{p\})$. Then we can use the long exact sequence associated to 
\eqref{blup1} to obtain information about $H^1(S', \T_{S'})$:  
\[ 
\begin{aligned} 
0 \lright &H^0(S', \T_{S'}) \lright H^0(S', \s^* \T_{S}) \lright 
H^0(\EE, \OO _{\EE} (1)) \lright \\ 
\lright &H^1(S', \T_{S'}) \lright H^1(S', \s^* \T_{S}) \lright H^1(\EE, \OO 
_{\EE} (1)) = 0. 
\end{aligned} 
\] 
 
Using Lemma \ref{P:blup}, we obtain: 
\begin{equation}\label{E:blup} 
\begin{aligned} 
H^0(S', \T_{S'})&=\ker(H^0(S, \T_{S})  \lright H^0(\EE, \OO _{\EE} 
(1))) \text{ and}\\  
h^1(S', \T_{S'}) &= 
\dim (\coker(H^0(S, \T_{S})  \lright H^0(\EE, \OO _{\EE} (1)))) 
+h^1(S, \T_{S})\\ 
&= h^0(S', \T_{S'})-h^0(S, \T_{S})+2+h^1(S, \T_{S}). 
\end{aligned} 
\end{equation}  
  
Note that if $S$ is obtained by successive blow-ups of (possibly 
infinitely near) points on $\PP^2$, then we can regard $H^0(S, \T_{S})$ as a subspace of $H^0 (\PP ^2, \T _{\PP ^2})$.  
 
The following theorem describes the relation between the first 
cohomology group of the tangent bundle of $V(\Sigma)$ and that of $V(\Sigma \cup \{p\})$. 
\begin{theorem}\label{T:DD} 
Let $\Sigma$ be a set of (possibly infinitely near) points of $\PP^2$ in almost general 
position. Suppose $|\Sigma| \leq 3$. Let $\sigma \colon S=V( \Sigma) \lright 
\PP^2$ be the blow-up of  center $\Sigma$. Let $p \in S$ be a 
point 
such that $\Sigma' = \Sigma \cup \{p\}$ is in almost general 
position and let $\s' \colon S' =V( \Sigma')\lright S$ be the blow-up of $p$ on
$S$. Denote by $\tilde \EE$ the union of all curves with negative self-intersection on $S$. Then we have the following:  
\begin{enumerate}  
   \item {\bf Type 1:} If $p \notin \tilde \EE$, then we have $h^1(S', \T_{S'}) = h^1(S, \T_{S})$.\label{type1}  
   \item {\bf Type 2:} If $p$ is contained in a single $(-1)$-curve, 
   then we have $h^1(S', \T_{S'}) = h^1(S, \T_{S}) + 1$.  
   \item {\bf Type 3:} If $p$ is the intersection point of two $(-1)$-curves, then $h^1(S', \T_{S'}) = h^1(S, \T_{S}) + 2$.  
\end{enumerate}  
\end{theorem} 

\begin{corollary}\label{C:DD} 
If $D$ is a normal rational del Pezzo surface of degree $d \geq 5$, 
with $r$ singularities, of type $A_{\lambda_1}$, $A_{\lambda_2}$, ..., 
 $A_{\lambda_r}$ respectively, and $\pi \colon \DD \lright D$ is its 
minimal resolution, then we have $\dim H^1 (\tilde D,\T_\DD)=\Sigma_{i=1}^{r} \lambda_i$. 
\end{corollary} 
\begin{proof} 
Note that $h^1(S', \T_{S'}) - h^1(S, \T_{S})$ counts the number of
$(-2)$-curves that appear on $S'$ after the blow-up $\s$, thus the corollary follows.
\end{proof} 
\begin{remark}  
For $d \leq 4$ (i.e. when $\DD$ is obtained by blowing up 5 or more
points) the Theorem is in general not true. For example, if we blow up 5
points in general position on $\PP^2$, the blow-up $\DD$ doesn't
contain (--2)-curves,
but $\dim H^1 (\tilde D,\T_\DD)=2$.
\end{remark}

\begin{proof}[Proof of Theorem \ref{T:DD}] 
Based on Remark \ref{R:config},  if $\Sigma$ and $\widetilde\Sigma$ (with 
$|\Sigma| =|\widetilde\Sigma| \leq 4$) have the same configuration
(i.e. if they can be transformed into each other by a projective
automorphism of $\PP^2$), then 
$V(\Sigma) \simeq V(\widetilde\Sigma)$. Now, if $p_1 \in V(\Sigma)$ and 
$p_2 \in V(\widetilde\Sigma)$ are points of the same type (1,2 or 3), then 
$\Sigma \cup \{p_1\}$ and $\widetilde\Sigma \cup \{p_2\}$ will again have 
the same configuration, and hence $V(\Sigma \cup \{p_1\}) \simeq 
V(\widetilde\Sigma \cup \{p_2\})$. Therefore we have to prove the theorem  
for only one representative $\Sigma$ for any configuration of points,
and we can use those described in \S \ref{delp}.  
 
We prove the theorem by blowing up one point at a time and explicitly 
computing the cohomologies involved. 
  
\chunk{\bf Explicit computation of $H^0 (\PP ^2, \T _{\PP ^2})$.}  
We know that $\dim H^0 (\PP ^2, \T _{\PP ^2})=~8$ and 
$\dim H^1 (\PP ^2, \T _{\PP ^2})=0$. Here we compute a 
basis for $H^0 (\PP ^2, \T _{\PP ^2})$ in local coordinates.  
  
Fix the homogeneous coordinates $[x_0 \colon x_1 \colon x_2]$ on 
$\PP^2$. Then on the affine open $U_0=\{x_0 \neq 0\}$ we have local 
coordinates $x: = \frac{x_1}{x_0}$ and $y := 
\frac{x_2}{x_0}$. Around $p=[1:0:0]$, $\T _{\PP ^2}$ is generated by 
the vectors $\frac {\partial}{\partial x}$, $\frac {\partial}{\partial 
y}$; more precisely, $\T _{\PP ^2}|_{U_0}=\C[x,y]\frac 
{\partial}{\partial x}+\C[x,y]\frac {\partial}{\partial y} \simeq \C[x,y]^2$. 
  
\begin{claim}  
With the above notations, $H^0 (\PP ^2, \T _{\PP ^2})$ has a basis given by 
 $(1,0)$, $(x,0)$, $(y,0)$, $(0,1)$, $(0,x)$, $(0,y)$, $(x^2,xy)$, $(xy,y^2)$ on 
 $U_0$. 
\end{claim}  
  
\begin{proof}  
We have the dual of the Euler sequence,  
\[  
0 \longrightarrow \OO _{\PP ^2} \longrightarrow \OO _{\PP ^2} (1)^3 \longrightarrow \T _{\PP ^2} \longrightarrow 0 ,  
\]  
where $\OO _{\PP ^2} (1)^3 \longrightarrow \T _{\PP ^2}$ is (locally) given by   
\[  
(s_0, s_1, s_2) \mapsto {\frac {s_1x_0 -s_0 x_1}{x_0 ^2}}{\frac {\partial}{\partial x}} + {\frac {s_2x_0 -s_0 x_2}{x_0 ^2}}{\frac {\partial}{\partial y}}.  
\]  
Writing out the generators of $H^0 (\PP ^2, \OO _{\PP ^2} (1)^3)$, the claim follows.  
\end{proof}  
Using this explicit description of $H^0 (\PP ^2, \T_{\PP^2})$, the computations are straightforward using the descriptions of the
blow-ups from \S \ref{delp}. We illustrate it in the case of
blowing up 1 point, respectively 2 infinitely near points on $\PP^2$, and leave the rest of
the computations to be  carried out by the reader.

In the following, all blow-ups will be of 
(possibly infinitely near) points on $U_0=\{x_0 \neq 0\} 
\subset \PP^2$ (as described in \S \ref{delp}). We will always denote by $\EE$ the exceptional locus 
of the last blow-up $S' \lright S$, and by $\mathcal M$ the map $\T _{S'}  \longrightarrow \sigma ^* \T_{S}$.

\chunk{\bf Blowing up a point.}\label{1}  
Let $\s_1 \colon S_1 \lright \PP^2$ be the blow-up of the point 
$(x,y)=(0,0)$. 
Based on \eqref{E:blup}, in order to find $\dim H^1(S_1, \T_{S_1})$, 
we first need to describe the kernel of the map $H^0(\PP^2, \T_{\PP^2}) \simeq H^0(S_1, \s_1^* \T_{\PP^2}) \lright H^0(\EE, \OO _{\EE} (1))$. Locally, near $\EE$, $S_1$ is covered by two affine opens:   
\[  
\begin{aligned}  
V_0&=\Spec \C[x,s], \text{where }x=x, y=xs\\  
V_1&=\Spec \C[y,t], \text{where }x=yt, y=y  
\end{aligned}
\]
The map $\mathcal M : \T_{S_1} \lright \sigma_1^* \T_{\PP^2}$ is given by the following matrices:  
\begin{equation*} 
\begin{aligned}  
\text{On } V_0  \/&:  \/ \mathcal M_0 =\begin{pmatrix}   
1 &0\\  
s &x  
\end{pmatrix}\text{ with basis elements $\frac {\partial}{\partial x}$, $\frac {\partial}{\partial s}$
for $\T_{S_1}|_{V_0}$ and $\frac {\partial}{\partial x}$, $\frac
{\partial}{\partial y}$ for $\s_1^* \T_{\PP^2}{|_{V_0}}$.}\\  
\text{On } V_1 \/&: \, \/ \mathcal M_1 =\begin{pmatrix}  
t &y\\  
1 &0  
\end{pmatrix}\text{ with basis elements $\frac {\partial}{\partial y}$, $\frac {\partial}{\partial t}$ for $\T_{S_1}|_{V_1}$.}  
\end{aligned}  
\end{equation*}

We have $\ker (H^0(S_1, \s_1^* \T_{\PP^2}) \rright H^0(\EE, \OO _{\EE} 
(1))) = \mathcal M \!(H^0(S_1, \T_{S_1}))\! \subset \!H^0(S_1, \s_1^* 
\T_{\PP^2}) \simeq H^0(\PP^2, \T_{\PP^2}).$  
 
To find this kernel, it is sufficient to find the elements of $H^0(\PP^2, 
\T_{\PP^2})$ that are in the image of $\mathcal M_0$. Indeed, if we 
have a global section $s \in H^0(\PP^2, \T_{\PP^2}) \simeq H^0(S_1, \s_1^* \T_{\PP^2})$  
in the kernel of $H^0(V_0, \s_1^* \T_{\PP^2}) \lright H^0(\EE 
\cap {V_0}, \OO _{\EE}(1))$, then $s$ is also  
in the kernel of the map $H^0(S_1, \s_1^* \T_{\PP^2}) \lright H^0(\EE, \OO _{\EE} 
(1))$, because $\OO _{\EE}(1)$ is a locally free sheaf on $\EE$, and 
hence if a global section vanishes on an open subset, then it vanishes 
on $\EE$. 
 
Locally on $V_0$ we have $\C[x,s]^2 \xrightarrow{\mathcal M_0}
\C[x,s]^2,$ and any element of $H^0\!(S_1, \s_1^* 
\T_{\PP^2}) \!\simeq H^0\!(\PP^2, \T_{\PP^2})$ is of form $\begin{pmatrix}  
a_1 +a_2 x + a_3 y +a_7 x^2 + a_8 xy\\  
a_4 +a_5 x + a_6 y +a_7 xy + a_8 y^2  
\end{pmatrix}$ where $y=xs$, and $a_i \in \C$.

We need the following lemma from algebra:  
\begin{lemma}\label{algebra}  
Let $A$ be an integral domain and $M \colon A^2 \lright A^2$ an $A$-linear map, given by the matrix $M =\begin{pmatrix}  
a &b\\  
c &d  
\end{pmatrix}  
$. Then an element $\begin{pmatrix}  
f\\  
g  
\end{pmatrix} \in A^2$ is in the image of $M$ if and only if $\dd=det(M)$ divides both $ag-cf$ and $df-bg$. 
\end{lemma}  
\begin{proof}  
The proof is obvious, using Cramer's rule in the quotient ring of $A$.  
\end{proof}

Suppose an element $\begin{pmatrix}  
a_1 +a_2 x + a_3 y +a_7 x^2 + a_8 xy\\  
a_4 +a_5 x + a_6 y +a_7 xy + a_8 y^2  
\end{pmatrix}$ of $H^0(\PP^2, \T_{\PP^2})|_{V_0}$ is in the image of 
$\mathcal M_0$. Then, by Lemma \ref{algebra}, $(a_4 +a_5 x + a_6 y 
+a_7 xy + a_8 y^2)-s(a_1 +a_2 x + a_3 y +a_7 x^2 + a_8 xy)$ is 
divisible by $x$, and this is equivalent, using the condition $y=xs$, 
to $a_4-sa_1=0$. This implies $a_4=a_1=0$. Hence  
\[  
 H^0(S_1,\T_{S_1}) \simeq \Big\{ \begin{pmatrix}  
a_1 +a_2 x + a_3 y +a_7 x^2 + a_8 xy\\  
a_4 +a_5 x + a_6 y +a_7 xy + a_8 y^2  
\end{pmatrix} \in H^0(\PP^2, \T_{\PP^2}) \colon a_4=a_1=0\Big\}  
\]  
has dimension 6, and therefore, by \eqref{E:blup}, $h^1(S_1,\T_{S_1})=h^1(\PP^2, \T_{\PP^2})=0$.

\chunk{\bf Blowing up two infinitely near points of $\PP^2$.}\label{22} 
Let $\s_7 \colon S_7 \lright S_1$ be the  blow-up of the point $(x,s)=(0,0) \in S_1$ (corresponding to the direction given by the line $\{y=0\}$ in $\PP^2$). Then near the exceptional locus $\EE$ of $\s_7$, $S_7$ is covered by the following affine opens:  
\[  
\begin{aligned}  
V_0&=\Spec \C[x,u], \text{where }x=x, y=xs=x^2u\\  
V_1&=\Spec \C[s,v], \text{where }x=sv, y=xs=s^2v\\  
V_2&=\Spec \C[y,t], \text{where }x=yt, y=y.  
\end{aligned}  
\]  
The exceptional locus $\EE$ is covered by $V_0$ and $V_1$. We compute on $V_0$. The corresponding matrix is 
\[  
\mathcal M_0 =\begin{pmatrix}  
1 &0\\  
2xu &x^2  
\end{pmatrix} \colon \C[x,u]^2 \lright \C[x,u]^2.  
\]  
Similar calculations as in Case \ref{1} lead to   
\[  
 H^0(S_7,\T_{S_7}) \simeq \Big\{ \begin{pmatrix}  
a_2 x + a_3 y +a_7 x^2 + a_8 xy\\  
a_5 x + a_6 y +a_7 xy + a_8 y^2  
\end{pmatrix} \in H^0(S_1,\T_{S_1}) \colon a_5=0\Big\} . 
\]  
Therefore $\dim  H^0(S_7,\T_{S_7})=5$ and so
$h^1(S_7,\T_{S_7})=h^1(S_1,\T_{S_1})+1=1$.   

For the rest of the cases, similar computations -using the description
of configurations of at most 4 almost general points from \S
\ref{delp}- yield to the proof of the theorem.
\end{proof}

\section{Example of a $K$-trivial contraction}\label{examples} 
 
Let $\varphi \colon X \longrightarrow~Y$ be a $K$-trivial birational 
extremal contraction of a smooth projective 3-fold $X$, contracting a divisor $D \subset X$ to a point $q \in Y$. Suppose $D$ is a 
normal rational del Pezzo surface of degree $d \geq 5$.  
The result of analytic rigidity shows that in order to know the 
analytic structure of the contraction $\f$, it is sufficient to have one 
example for each possible exceptional divisor with the prescribed 
normal bundle.

In the following,  we construct an example  of 
embedding a normal rational (singular) del Pezzo surface $D$ of degree
7 into a smooth 
threefold $X$ with the prescribed normal bundle $\OO_D(K_D)$, and  hence, by Fujiki's  
contraction theorem (see Theorem \ref{T:Fujiki} below), obtain an analytic 
contraction of $D$ (i.e. a holomorphic map $\f \colon X \lright Y$ onto a 
normal complex space $Y$ that contracts $D$ to a point $q \in Y$). Similar constructions can be
carried out  for each possible exceptional 
divisor $D$ (any normal rational del Pezzo 
  surface of degree $d \geq 5$).

Let $D$ be a normal rational del Pezzo surface of degree $d=7$. 
 By the results
in \S \ref{delp}, the minimal resolution $\DD$ of $D$ is isomorphic to
$V(\{p_1,p_2\})$, where $p_2 \lright p_1$, $p_1=[1:0:0] \in \PP^2$
and $p_2=[1:0] \in \PP(T_{p_1} \PP^2 )$.

First  we construct a family $\X \lright \mathbb A^1$ 
such that $\X$ is nonsingular, the central fiber $\X_0$ is isomorphic
to $D$, and the general fibers of the 
family are nonsingular del Pezzo surfaces of degree $d=7=\deg D$.  We will construct the family
$\X \subset \PP^7 \times \mathbb A^1$ as the closure of a
family of blow-ups of 2 distinct points on $\PP ^2$. 

For this, consider the curve $\mathcal C =\{x_1^2-t(x_1+x_0)=0, x_2=0\}
\subset \PP^2 \times \mathbb A^1$. Over any $t \neq 0$, $\mathcal C$ has two
points, while over $t=0$ it has a double point. The curve $\mathcal C$
defines a map $\Phi =\{\varphi _t \}_t: \PP ^2 \times \mathbb A^1\cdots \rightarrow \PP ^7\times \mathbb A^1$ with
base-locus $\mathcal C$, where
$\varphi _t: \PP ^2 \cdots \rightarrow \PP ^7$ is given by 
\[
[x_0:x_1:x_2] \longmapsto 
[x_1^3-tx_0^2x_1-t^2x_0^3-t^2x_0^2x_1:x_2^3:x_0^2x_2:x_0x_1^2-tx_0^3-tx_0^2x_1:x_0x_2^2:x_1^2x_2:x_1x_2^2:x_0x_1x_2].
\]

Let $\X$ be the
blowing up of $\PP^2 \times \mathbb A^1$ with center $\mathcal C$. Then
$\f_{|-K_{\X/\mathbb A^1}|}$ will define an embedding of $\X$ into $\PP
^7\times \mathbb A^1$. The ideal $\I$ defining $\X$ in this embedding
is
\begin{center} 
$\I=(ge-bh, -hf+gd+th^2+{ tch}, -he+gc, f^2-ga-t^2h^2 -{th^2-t^2ch}, ef-hg,$\\ 
 $fd-ah+thd, fc-h^2, -g^2+bf, de-h^2+{{tc^2}}+tch, ea-hf+t^2ch+{t^2c^2+tch},$\\$bd-hg+teh+{tce}, -e^2+cb, ca-dh-tcd, ba-fg+t^2eh+{t^2ce+tcg}).$ 
\end{center}

For each $t$, $\I$ defines a surface $D_t$ in $\PP ^7$. It can be verified by direct
computation that the total space $\X$ of the 
family $\{D_t\}_{t \in \mathbb A^1}$ is nonsingular and that $D_t$ is
nonsingular (actually a smooth del Pezzo surface of degree 7) except
in the case of
$D_0=D$ and $D_{-4}$, which are singular del Pezzo surfaces. In fact,
for $t \neq 0,-4$,  the surface  $D_t$ is the image (via the
anticanonical embedding) in $\PP ^7$ of the blow up of the
points  $[1:\lambda_{1}:0]$ and $[1:\lambda_{2}:0]$ on $\PP ^2$, where
$\lambda _{1}$ and $\lambda _{2}$ are the roots of
the equation $u^2=t(u+1)$.

Since $D$ is a fiber of the family $\X$, 
its normal bundle $\N_{D/\X} \simeq \OO_D(D)$ in $\X$ is numerically trivial.  

In order to make the normal bundle isomorphic to $\OO_D(K_D)$, as
required for a $K$-trivial contraction, we proceed as follows: 
  
Let $C \in |-K_D|$ be a general member avoiding the singular point of 
$D$ (it exists, because $|-K_D|$ is very ample). Let $X$ be the 
blow-up of $C$ on $\X$, $E$ the exceptional locus, $D'$ the strict 
transform of $D$, and $C' := E|_D' \simeq C$. Because $K_{\X}|_D$ and 
$K_D$ are numerically equivalent, it follows that $K_{X}$ is 
numerically trivial on $D'$. Therefore $\N_{D'/X}$ is numerically equivalent 
to $\OO_{D'}(K_{D'})$, and so, as in \S \ref{known}, $\N_{D'/X} \simeq 
\OO_{D'}(K_{D'})$. Therefore we succeeded in embedding the del Pezzo surface $D 
\simeq D'$ into a nonsingular three-fold $X$ such that its normal 
bundle is isomorphic to the canonical sheaf of $D$. We can now apply the 
following theorem (\cite[Th.2.]{Fujiki}): 
  
\begin{theorem}[Fujiki's Contraction Theorem]\label{T:Fujiki} 
Let $X$ be a complex space, $A \subset X$ an effective Cartier divisor, $B$ another complex space, and $f:A \lright B$ a surjective holomorphic map. Assume that:  
1. The conormal bundle $\N_{A/X}^*$ is f-positive, and that  
2. $R^1 f_* ({\N_{A/X}^*}^{\nu})=0$, for all $\nu >0$.  
Then there exists a modification $\psi :X \lright Y$ with $\psi |_A=f$. Moreover, $\psi _*(\mathcal L) \simeq \OO_Y$, where the coherent sheaf $\mathcal L$ is defined by   
\[  
0 \lright \mathcal L \lright \OO_{\X} \lright \OO_A/\im(f^*  \OO_B \lright \OO_A) \lright 0.  
\]  
\end{theorem}   
  
In our case $A:=D'\simeq D$, $B$ is a point (hence $\mathcal L \simeq 
\OO_{X}$) and $R^1 f_* ({\N_{A/X}^*}^\nu) \simeq H^1(D,\OO_D(\nu K_{D}))=0$ for all 
$\nu$ (see Proposition \ref{P:cohom}). 
  
Therefore we have an analytic modification $\psi : X \lright Y$
contracting $D$ to a point $q \in Y$, such that $\psi |_{X-D} : X-D \simeq Y-q$, and $\psi _* \OO_{X} =\OO_Y$.  
 
\begin{remark} 
A similar result is true in the category of algebraic spaces, based on 
a contraction theorem due to Artin (\cite[Cor. (6.10)]{Artin2}). 
\end{remark} 
However, Fujiki's theorem doesn't guarantee a contraction in the 
algebraic category (i.e. the existence of a morphism $\f \colon X 
\lright Y$ of algebraic varieties contracting $D$ to a point). At 
the present we do not know how to construct in general algebraic morphisms 
that contract a singular del Pezzo surface embedded into a smooth 
3-fold with the prescribed normal bundle to a point, i.e. how to 
obtain a $K$-trivial morphism $\f \colon X \lright Y$ onto a 
normal projective variety $Y$ that contracts $D$ to a point $q \in 
Y$. In the case when $D$ is nonsingular, the contraction of the 
zero-section of the total space of the normal bundle $\N_{D/X}$ 
provides such an example. Although this gives a $K$-trivial contraction, 
it is not extremal. Namikawa constructs an 
example of a $K$-trivial {\it extremal} contraction (\cite[Example 
1]{Nam1}) with exceptional divisor $D$ a nonsingular del Pezzo 
surface of degree 6 (i.e. a smooth cubic).

\begin{acknowledgement}
{The author would like to thank K. Matsuki for his advice and
encouragement, and also M. Gross and J. Koll\'ar for helpful
comments.}
\end{acknowledgement}

\end{document}